\def\isarxiv{1}

\def\paperTitle{Differential Privacy for Euclidean Jordan Algebra with Applications to Private Symmetric Cone Programming}

\def\paperAuthor{
Zhao Song\thanks{\texttt{magic.linuxkde@gmail.com}. University of California, Berkeley.}
\and
Jianfei Xue\thanks{\texttt{jx898@nyu.edu}. New York University.} \and
Lichen Zhang\thanks{\texttt{lichenz@mit.edu}. Massachusetts Institute of Technology. Supported by a Mathworks Fellowship and a Simons Dissertation Fellowship in Mathematics.}
}

\ifdefined\isarxiv
\documentclass[11pt]{article}
\usepackage[numbers]{natbib}
\else
\documentclass{article}
\usepackage[final]{neurips_2025}
\fi

\ifdefined\isarxiv

\usepackage{amsmath}
\usepackage{amsthm}
\usepackage{amssymb}
\usepackage{algorithm}
\usepackage{subfig}
\usepackage{algpseudocode}
\usepackage{graphicx}
\usepackage{grffile}
\usepackage{wrapfig,epsfig}
\usepackage{url}
\usepackage{xcolor}
\usepackage{epstopdf}

\usepackage{bbm}
\usepackage{dsfont}

\else 

\usepackage[utf8]{inputenc} 
\usepackage[T1]{fontenc}    
\usepackage{hyperref}       
\usepackage{url}            
\usepackage{booktabs}       
\usepackage{amsfonts}       
\usepackage{nicefrac}       
\usepackage{microtype}      
\usepackage{xcolor}         

\fi

\ifdefined\isarxiv

\else
\usepackage{amsmath}
\usepackage{amsthm}
\usepackage{amssymb}
\usepackage{algorithm}
\usepackage{subfig}
\usepackage{algpseudocode}
\usepackage{graphicx}
\usepackage{grffile}
\usepackage{wrapfig,epsfig}
\usepackage{epstopdf}
\usepackage{bbm}
\usepackage{dsfont}
\fi
 
\allowdisplaybreaks

\ifdefined\isarxiv

\usepackage{tikz}
\usepackage{hyperref}  
\hypersetup{colorlinks=true,citecolor=blue,linkcolor=blue} 
\usetikzlibrary{arrows}
\usepackage[margin=1in]{geometry}

\else
\definecolor{mydarkblue}{rgb}{0,0.08,0.45}
\hypersetup{colorlinks=true, citecolor=mydarkblue,linkcolor=mydarkblue}
\fi
 
\graphicspath{{./figs/}}

\theoremstyle{plain}
\newtheorem{theorem}{Theorem}[section]
\newtheorem{lemma}[theorem]{Lemma}
\newtheorem{definition}[theorem]{Definition}

\newtheorem{corollary}[theorem]{Corollary}

\newtheorem{remark}[theorem]{Remark}

\newcommand{\wh}{\widehat}
\newcommand{\wt}{\widetilde}
\newcommand{\ov}{\overline}
\newcommand{\N}{\mathcal{N}}
\newcommand{\R}{\mathbb{R}}

\newcommand{\K}{\mathcal{K}}
\newcommand{\J}{\mathcal{J}}

\DeclareMathOperator{\OPT}{OPT}

\DeclareMathOperator{\poly}{poly}

\DeclareMathOperator{\rank}{rank}

\DeclareMathOperator{\Tr}{Tr}

\newenvironment{CompactItemize}{
\begin{list}{\tiny$\bullet$}{%
\setlength{\leftmargin}{10pt}
\setlength{\itemindent}{0pt}
\setlength{\topsep}{-1pt}
\setlength{\itemsep}{0pt}
}}
{\end{list}}

\begin{document}

\ifdefined\isarxiv

\date{}
\title{\paperTitle}
\author{\paperAuthor}

\else

\title{\paperTitle}

\author{%
  Zhao Song \\
  University of California, Berkeley \\
  \texttt{magic.linuxkde@gmail.com} \\
  \And 
  Jianfei Xue \\
  New York University \\
  \texttt{jx898@nyu.edu} \\
  \And
  Lichen Zhang \\
  MIT CSAIL \\
  \texttt{lichenz@csail.mit.edu} \\
}

\maketitle

\fi

\ifdefined\isarxiv
\begin{titlepage}
  \maketitle
  \begin{abstract}
    In this paper, we study differentially private mechanisms for functions whose outputs lie in a Euclidean Jordan algebra. Euclidean Jordan algebras capture many important mathematical structures and form the foundation of linear programming, second-order cone programming, and semidefinite programming. Our main contribution is a generic Gaussian mechanism for such functions, with sensitivity measured in $\ell_2$, $\ell_1$, and $\ell_\infty$ norms. Notably, this framework includes the important case where the function outputs are symmetric matrices, and sensitivity is measured in the Frobenius, nuclear, or spectral norm. We further derive private algorithms for solving symmetric cone programs under various settings, using a combination of the multiplicative weights update method and our generic Gaussian mechanism. As an application, we present differentially private algorithms for semidefinite programming, resolving a major open question posed by [Hsu, Roth, Roughgarden, and Ullman, ICALP 2014].

  \end{abstract}
  \thispagestyle{empty}
\end{titlepage}

\newpage

\else

\begin{abstract}
In this paper, we study differentially private mechanisms for functions whose outputs lie in a Euclidean Jordan algebra. Euclidean Jordan algebras capture many important mathematical structures and form the foundation of linear programming, second-order cone programming, and semidefinite programming. Our main contribution is a generic Gaussian mechanism for such functions, with sensitivity measured in $\ell_2$, $\ell_1$, and $\ell_\infty$ norms. Notably, this framework includes the important case where the function outputs are symmetric matrices, and sensitivity is measured in the Frobenius, nuclear, or spectral norm. We further derive private algorithms for solving symmetric cone programs under various settings, using a combination of the multiplicative weights update method and our generic Gaussian mechanism. As an application, we present differentially private algorithms for semidefinite programming, resolving a major open question posed by [Hsu, Roth, Roughgarden, and Ullman, ICALP 2014].

\end{abstract}

\fi


\section{Introduction}

As modern machine learning leverages increasingly large models and ever-growing datasets, protecting the privacy of data has become a paramount concern. Differential privacy~\cite{dmns06,dkmmn06} has emerged as the gold standard for data privacy, owing to its simplicity, practicality, and strong theoretical guarantees. Designing machine learning algorithms with differential privacy (DP) is particularly valuable, as it enables protection of both user data and model parameters~\cite{bik+17,gkn17,mrtz18,tlc+20,atmr21,nbd22,fhl+24}. However, most differentially private algorithms for machine learning are developed in a case-by-case manner — tailored specifically to individual problems to ensure good privacy and utility trade-offs. More generic approaches, such as making stochastic gradient descent (SGD) differentially private~\cite{scs13,bst14,acg+16}, have gained significant traction. Since SGD is a core subroutine in many training algorithms, this line of work provides automatic privacy guarantees across a wide range of models. Following this direction, we develop differentially private algorithms for solving \emph{symmetric cone programs} (SCP), a broad class of convex optimization problems that includes linear programs (LP), second-order cone programs (SOCP), and semidefinite programs (SDP). These convex formulations arise in many machine learning applications, including support vector machines~\cite{bgv92,cv95,mmrts01,j06,cl11,sssc11,gsz25}, matrix completion~\cite{ct10,r11,cr12,jns13,zwl15,glz17,kllst23,gsyz24,syyz25}, robust mean and covariance estimation~\cite{cdg19,cdgw19}, experimental design~\cite{vg19,azlsw20}, and sparse PCA~\cite{dejl04,zht06,dbe08,zde10,vclr13}.

Symmetric cone programming (SCP) has been extensively studied from an algorithmic perspective, with both first-order methods~\cite{clpv23,zvtl24} and second-order methods~\cite{f97_quadratic,f97_linear,sa03,p23} developed, as they provide valuable insights for solving semidefinite programs. However, algorithms for SCP with differential privacy (DP) guarantees have received little attention, largely due to the abstract nature of the problem class. In fact, to the best of our knowledge, the only known DP algorithms in this family are for linear programs~\cite{hrru14,dfv+20,muvv21,bbdh24,bbdh25}, and designing DP algorithms for SDPs was posed as a major open question in~\cite{hrru14}. In this paper, we take a significant step forward by designing DP mechanisms not only for SDPs, but for the broader class of SCPs — a much richer and more general family of convex optimization problems. Our main contributions are two-fold: (1) We begin by examining Euclidean Jordan algebras (EJAs), a class of real inner product spaces that includes $\mathbb{R}^k$ and the space of real or complex symmetric matrices of size $r \times r$. We develop a generic Gaussian mechanism over EJAs that provides $(\epsilon,\delta)$-differential privacy guarantees when sensitivity is measured in the $\ell_1$, $\ell_2$, or $\ell_\infty$ norms. In the case of symmetric matrices, this corresponds to nuclear, Frobenius, and spectral norms, respectively. This stands in contrast to the well-known matrix mechanism~\cite{lhr+10,ht09}, which aggregates multiple linear queries into a matrix and measures sensitivity column-wise in vector norms. Our framework instead captures the geometric structure of matrices by modeling perturbations to their spectrum. (2) We then design a suite of private algorithms for solving SCPs under various neighboring data settings. Specifically, given two neighboring databases $D$ and $D'$, we consider the following cases: (i) $D'$ has one additional linear constraint compared to $D$ (high sensitivity constraint privacy), (ii) the scalar vectors $b$ differ by at most a sensitivity bound $\Delta_\infty$ in $\ell_\infty$ norm (low sensitivity scalar privacy), (iii) the constraint sets differ by at most $\Delta_\infty$ in $\ell_\infty$ norm (low sensitivity constraint privacy), and (iv) the objective elements differ by at most $\Delta_\infty$ in $\ell_\infty$ norm (low sensitivity objective privacy). For the high sensitivity constraint setting, we show that the multiplicative weights update algorithm over constraints from~\cite{hru13,hrru14} naturally extends to SCP when the optimal value $\OPT$ is given as a public input. For covering SDPs over $r\times r$ matrices, it outputs a solution of trace $\OPT$ such that, with high probability, fewer than $s$ constraints are violated by more than an additive $\alpha$, where $s$ grows linearly in $r$ and, up to logarithmic factors, as $\OPT^2/(\epsilon\alpha^2)$. For cases (ii) and (iii), we develop a novel variant of the multiplicative weights algorithm that operates over the primal variable space, guided by an approximate oracle that identifies the most violated constraint. In case (iii), this is combined with the Gaussian mechanism we previously designed for the EJA structure. In case (iv), we perturb the objective element with this Gaussian mechanism and solve the perturbed program with an added unit $\ell_2$ norm constraint, assuming an optimal solution of unit $\ell_2$ norm. As a byproduct, we design a multiplicative weights algorithm that could utilize a \emph{noisy oracle} to for solving SCPs, which we believe could lead to applications beyond differential privacy.

\paragraph{Roadmap.} In Section~\ref{sec:relat}, we survey related work on both differential privacy and symmetric cone programming. In Section~\ref{sec:preli}, we provide preliminary background on differential privacy and Euclidean Jordan algebras. Section~\ref{sec:dp_eja} introduces a generic Gaussian mechanism for Euclidean Jordan algebras and establishes privacy guarantees under $\ell_1$, $\ell_2$, and $\ell_\infty$ sensitivities. In Section~\ref{sec:dp_scp}, we present algorithms for private symmetric cone programming, covering both high sensitivity and low sensitivity constraint privacy settings. For the high-sensitivity case, we give algorithms for solving covering semidefinite programs and covering symmetric cone programs. Section~\ref{sec:concl} concludes with a discussion of open problems and limitations of our work. Appendix~\ref{sec:more_preli} provides additional preliminaries, particularly on Euclidean Jordan algebras and regret bounds for symmetric cone multiplicative weights updates. Appendix~\ref{sec:high_sensitiviy} presents the full details, analysis, and applications for the high sensitivity constraint privacy setting. In Appendix~\ref{sec:low_sensitivity}, we describe algorithms for the low sensitivity setting, including scalar privacy, constraint privacy, and objective privacy.

\section{Related Work}
\label{sec:relat}
\noindent{\bf Differential Privacy.} Since its introduction in~\cite{dmns06}, differential privacy has become the standard notion for providing rigorous privacy guarantees in algorithm design. It has found widespread applications across general machine learning~\cite{cm08,wm10,je19,tf20}, deep neural networks~\cite{acg+16,bps19}, computer vision~\cite{zycw20,lwaff21,tzxl19}, natural language processing~\cite{ydw+21,wk18}, federated learning~\cite{bik+17,syy+22,swyz23} and adaptive data structures~\cite{hkm+22,bkm+21,acss23,csw+23,syyz23_dp,ffl+25}. The line of work most relevant to ours concerns differentially private algorithms for convex programming, particularly linear programming. The work of~\cite{hru13} was among the first to study private algorithms for zero-sum games, a special case of linear programming. Subsequently,~\cite{hrru14} systematically examined private linear programs under various notions of sensitivity, including constraint, scalar, and objective perturbations. Their algorithms typically produce solutions that may violate constraints slightly due to privacy noise. To address this,~\cite{muvv21} considered the case where only the scalar vector is private, and showed that constraint satisfaction can be maintained by explicitly perturbing the scalar vector. An alternative approach was proposed in~\cite{dfv+20}, which reformulates the private problem as a stochastic chance-constrained program, whose solution satisfies the original constraints with high probability. Building on these ideas,~\cite{bbdh24} extended the framework of~\cite{muvv21} to handle private constraints, and~\cite{bbdh25} further demonstrated that it is possible to privatize all components of a linear program while still satisfying the constraints. Our work is most closely related to~\cite{hrru14}, as the other approaches heavily rely on the specific structure of linear programs and extensively use the Laplace mechanism. In contrast, we show that extending these ideas to symmetric cone programming poses significantly greater challenges, both technically and algorithmically.

\vspace{2mm}
\noindent{\bf Euclidean Jordan Algebras and Symmetric Cone Programming.} 
Euclidean Jordan algebras form the algebraic foundation of symmetric cone programming. It is well known that any EJA can be decomposed into a direct sum of simple Jordan algebras, and up to isomorphism, there are exactly five types of simple EJAs: real symmetric matrices, complex Hermitian matrices, quaternionic Hermitian matrices, spin factors, and the Albert algebra~\cite{fk94}. Among these, real symmetric matrices underlie semidefinite programming, while spin factors correspond to second-order cone programming. EJAs are commutative but non-associative algebras, and this structure provides key insights into the geometry and algorithmic design of semidefinite programming. As symmetric cone programming generalizes both SOCP and SDP, interior point methods have been extensively developed for this broader setting~\cite{f97_quadratic,f97_linear,sa03,v07}, with significant effort devoted to abstracting and extending classical SDP techniques. In particular,~\cite{p23} proposes a novel interior point method for symmetric cone programs using geodesic updates on the Riemannian manifold defined by the interior of the symmetric cone. Our work builds heavily on the multiplicative weights update framework introduced in~\cite{clpv23} for online linear optimization over symmetric cones, and extended in~\cite{zvtl24} to general symmetric cone programming. This framework is especially well suited for designing differentially private algorithms, as the oracle component required by multiplicative weights can be implemented via the exponential mechanism.

\paragraph{Subsequent work.} After the first version of this paper, Chu and Zhang~\cite{cz26} extended our high sensitivity results beyond covering programs and improved their dependence on the scale and the accuracy from polynomial to logarithmic. They study approximate feasibility of symmetric cone programs, including linear, second-order cone and semidefinite programs, under high sensitivity row privacy, where each record is one constraint, as in our high sensitivity constraint privacy. Their $(\epsilon,\delta)$-differentially private algorithm outputs a point in the cone that violates $\widetilde O(k^2/\epsilon)$ constraints by more than an additive error $\eta$, where $k$ is the dimension of the ambient space. This count depends only logarithmically on $UR/\eta$, where $R$ is the norm of a feasible point and $U$ is the size of the problem. In comparison, for covering SDPs we prove in Theorem~\ref{thm:private_covering_sdp} a bound that grows as $(\OPT/\alpha)^2$ but only linearly in $r$. Since $k=r(r+1)/2$ for $r\times r$ matrices, their bound with $\eta=\alpha$ is $\widetilde O(r^4/\epsilon)$, which is smaller than our $\widetilde O(r\cdot\OPT^2/(\epsilon\alpha^2))$ once $\OPT/\alpha$ exceeds $r^{1.5}$, up to logarithmic factors. They also prove that, even with additive error, a private algorithm must violate $\Omega(k/\epsilon)$ constraints on some linear programs.

\section{Preliminaries}
\label{sec:preli}

In this section, we provide some preliminary knowledge on Euclidean Jordan algebras, symmetric cone programming and differential privacy.
\subsection{Notations}
We use $\wt O(f)$ to denote $O(f\poly\log f)$. For two real vectors $x, y\in \R^k$, we use $\langle x,y\rangle$ to denote $x^\top y$. We use $\|x\|_1,\|x\|_2,\|x\|_\infty$ to denote the vector $\ell_1,\ell_2$ and $\ell_\infty$ norms. For two real symmetric matrices $X, Y\in \R^{r\times r}$, we use $\langle X, Y\rangle$ to denote $\Tr(X^\top Y)$, where $\Tr(\cdot)$ is the trace of the matrix. We use $\|X\|_{S_1},\|X\|_F$ and $\|X\|$ to denote the matrix nuclear, Frobenius and spectral norms. We use $X\succeq Y$ to denote $X-Y$ is a positive semidefinite matrix. We use $\N(\mu, \Sigma)$ to denote the multivariate Gaussian distribution with $\mu$ and covariance matrix $\Sigma$. 

\subsection{Euclidean Jordan Algebras and Symmetric Cone Programming}

We will exclusively work with Euclidean Jordan algebras, which, as we will see later, generalize many important spaces, including $\mathbb{R}^k$ and the set of all $r \times r$ real symmetric matrices.
\begin{definition}[Euclidean Jordan algebra (EJA)]
    An Euclidean Jordan algebra (EJA) is a finite-dimensional vector space $\J$ equipped with:
    \begin{CompactItemize}
        \item a bilinear
    product $\circ : \J \times \J \to \J$ that satisfies, for all $x, y \in \J$, $x \circ y = y \circ x$, and $x^2 \circ (x \circ y) = x \circ (x^2 \circ y)$;
    \item an inner product $\langle \cdot , \cdot \rangle :  \J \times \J \to \R$ that satisfies, for all $x, y, z\in \J$, $\langle x \circ y, z \rangle = \langle x , y \circ z \rangle$;
    \item an identity element $e$ that satisfies, for all $x\in \J$, $e \circ x = x \circ e = x$.
    \end{CompactItemize}
\end{definition}

As an example, if $\J$ is the set of all $r \times r$ real symmetric matrices, then the Jordan product is defined by $x \circ y = \frac{1}{2}(x y + y x)$, where $\cdot$ denotes standard matrix multiplication. Euclidean Jordan algebras induce a geometric structure known as a \emph{symmetric cone}, which can be characterized as follows:
\begin{definition}[Symmetric cone,~\cite{v07}]
    A symmetric cone is a closed convex cone $\K$ in a finite-dimensional inner product space $\J$ that satisfies the following properties:
    \begin{CompactItemize}
        \item $\K$ is self-dual, i.e., $\K^* := \{ y \in \J : \langle y, x \rangle \geq 0, \forall x \in \K \} = \K$.
        \item $\K$ is homogeneous, i.e., for any $u, v \in \mathrm{int}(\K)$, there exists an invertible linear transformation $L: \J \to \J$ such that $L(u) = v$ and $L(\K) = \K$.
    \end{CompactItemize}
    $\K$ can be characterized as the cone of squares of an EJA, specifically, there exists an EJA $\J$ with $\K=\{x^2: x\in \J\}$ where $x^2=x\circ x$ for $x\in \J$.
\end{definition}

There are two key parameters associated with a Euclidean Jordan algebra (EJA) $\J$: its \emph{dimension} and \emph{rank}.

\begin{definition}[Dimension and isomorphism]
\label{def:dim_iso}
Let $\J$ be an EJA. Then there exists a positive integer $k$ such that $\J$ is isomorphic to $\mathbb{R}^k$. This integer $k$ is called the \emph{dimension} of $\J$, denoted by $\dim(\J) = k$. Moreover, there exists a linear isomorphism $\phi : \J \rightarrow \mathbb{R}^k$ such that for all $x, y \in \J$, we have $\langle x, y \rangle = \langle \phi(x), \phi(y) \rangle$, i.e., $\phi$ is an isometry with respect to the inner product.
\end{definition}

As an example, consider the set of $r \times r$ real symmetric matrices. The dimension of this space is $\frac{r(r+1)}{2}$, since one can map these matrices to vectors of the same dimension. To define the rank, we first introduce a suitable spectral decomposition for any element $x \in \J$.

\begin{definition}[Jordan frame]\label{def:jordan_frame}
Let $q_1, \dots, q_r \in \J$ be elements satisfying idempotency ($q_i^2 = q_i$) and primitiveness ($q_i \neq 0$ and cannot be written as a sum of two nonzero idempotents). We say $\{q_1, \ldots, q_r\}$ is a Jordan frame if (1) $q_i \circ q_j = 0$ for $i \neq j$; (2) $\sum_{i=1}^r q_i = e$.
\end{definition}

The rank of an EJA is defined as the size of a Jordan frame used to form a \emph{spectral decomposition}.

\begin{definition}[Spectral decomposition]\label{lem:spectral_decomposition}
For any $x \in \J$, there exists a set of unique real numbers $\lambda_1, \ldots, \lambda_r$ and a Jordan frame $\{q_1, \ldots, q_r\}$ such that $x = \sum_{i=1}^r \lambda_i q_i$. The minimal such $r$ is called the rank of $\J$ (denoted by $\rank(\J) = r$). We use $\Tr(x) = \sum_{i=1}^r \lambda_i$ to denote the trace of the element $x$. We define the exponentiation function $\exp : \J \rightarrow \J$ as $\exp\left(\sum_{i=1}^r \lambda_i q_i\right) := \sum_{i=1}^r \exp(\lambda_i) q_i$.
\end{definition}

Returning to the example of the set of $r \times r$ real symmetric matrices, we observe that the spectral decomposition coincides with the standard matrix spectral decomposition. Therefore, the rank of this space is $r$. Note the quadratic gap between the rank and the dimension for this example. The final concept we introduce here is the trace-based inner product and the associated family of norms.

\begin{definition}[Inner product and norm]\label{def:inner_product}
Given $x, y \in \J$, we define their inner product as $\langle x, y \rangle = \Tr(x \circ y)$. For any $p \in [1, \infty)$, the $\ell_p$ norm of $x$ is defined as $\|x\|_p = \left(\sum_{i=1}^r |\lambda_i(x)|^p\right)^{1/p}$, and for $p = \infty$, the $\ell_\infty$ norm is defined as $\|x\|_\infty = \max_{i \in [r]} |\lambda_i(x)|$.
\end{definition}

Standard inequalities related to the inner product, such as Cauchy–Schwarz and H{\"o}lder's inequality, continue to hold in this setting; we defer their proofs to Appendix~\ref{sec:more_preli}. We are now ready to define symmetric cone programming.

\begin{definition}[Symmetric cone programming]\label{def:symmetric_cone_program}
Given a collection of elements $a_1, \ldots, a_m, c \in \J$ and $b \in \mathbb{R}^m$, the symmetric cone program (SCP) is defined as
\begin{align*}
    \max & ~\langle c, x \rangle \\
    \textnormal{s.t.} & ~\langle a_i, x \rangle \leq b_i, \quad \forall i \in [m], \\
    & ~ x \succeq_{\K} 0,
\end{align*}
where $\succeq_{\K}$ denotes the generalized cone inequality, i.e., $x \succeq_{\K} y$ means $x - y \in \K$.
\end{definition}

\subsection{Differential Privacy}
Differential privacy is a strong and rigorous notion of privacy, first introduced in~\cite{dmns06}. We state its formal definition below.

\begin{definition}[Differential privacy,~\cite{dmns06}]\label{def:differential_privacy}
A randomized mechanism ${\cal M} : \mathcal{D} \to \mathcal{R}$ is said to provide 
$(\epsilon,\delta)$-differential privacy if, for every pair of neighboring databases 
$D,D' \in \mathcal{D}$ (differing in exactly one record) and every measurable subset 
$S \subseteq \mathcal{R}$,
\begin{align*}
    \Pr[{\cal M}(D) \in S] \leq & ~ e^{\epsilon} \Pr[{\cal M}(D') \in S] + \delta .
\end{align*}
Here $\epsilon > 0$ and $0 \leq \delta < 1$ are privacy parameters. 
The special case $\delta = 0$ is referred to as $\epsilon$-differential privacy.
\end{definition}

We will make use of the standard Gaussian mechanism.

\begin{definition}[Gaussian mechanism~\cite{dkmmn06}]
\label{def:gaussian_mech}
Let $\epsilon, \delta \in (0,1)$. A mapping $f : \mathcal{D} \to \mathbb{R}^k$ has $\ell_2$-sensitivity $\Delta_2$ 
if for all neighboring databases $D,D'$ differing in one record, $\| f(D) - f(D') \|_2 \leq \Delta_2$. The Gaussian mechanism releases ${\cal M}(D) = f(D) + \nu$, where $\nu \sim \mathcal{N}(0, \sigma^2 I_k)$ and $\sigma = \frac{\Delta_2 \sqrt{2 \log(1.25/\delta)}}{\epsilon}$. This mechanism satisfies $(\epsilon,\delta)$-differential privacy.
\end{definition}

We note that~\cite{bw18} has improved the parameters of Gaussian mechanism, in particular the dependence on $\epsilon, \delta$ for $\sigma$. They further extend the analysis to work for all ranges of $\epsilon$. In this work, we focus on providing a preliminary privacy analysis of EJA, thus we adopt the more traditional and simple bound for Gaussian mechanisms.

To analyze the utility of the Gaussian mechanism, it is convenient to recall a standard tail bound for Gaussian random vectors.

\begin{lemma}[Laurent and Massart~\cite{lm00}]\label{lem:chi_square_tail}
Let $Z \sim \chi_k^2$ be a chi-squared random variable with $k$ degrees of freedom. Suppose each component has zero mean and variance $\sigma^2$. Then:
\begin{align*}
\Pr[ Z - k \sigma^2 \geq (2 \sqrt{k t} + 2t) \sigma^2 ] \leq & ~ \exp(-t), \\
\Pr[ k \sigma^2 - Z \geq 2 \sqrt{k t} \sigma^2 ] \leq & ~ \exp(-t).
\end{align*}
\end{lemma}

We also make use of the exponential mechanism~\cite{mt07}, a tool for achieving privacy when the output lies in a discrete or non-numeric domain. This mechanism relies on a \emph{quality score}, a real-valued function that evaluates the utility of pairing a database with a candidate output. Given a database, the exponential mechanism then selects, in a privacy-preserving way, an element whose quality score is close to the maximum achievable value.

\begin{definition}[Exponential mechanism,~\cite{mt07}]\label{def:exponential_mechanism}
For a database $D$ and range $\mathcal{R}$, the exponential mechanism 
chooses $r \in \mathcal{R}$ with probability proportional to $\exp \left( \frac{\epsilon}{2\Delta}  Q(r,D) \right)$. Here $Q : \mathcal{R} \times \mathcal{D} \to \mathbb{R}$ is a quality score 
function, and $\Delta$ denotes its sensitivity with respect to neighboring databases. This mechanism is $\epsilon$-differential privacy.
\end{definition}

The exponential mechanism satisfies the following utility guarantee.

\begin{lemma}[Accuracy guarantee for exponential mechanism,~\cite{kpru14}] \label{lem:acc_gua_exp_mec}
Given $\beta\in (0,1)$ and database $D$, let $\OPT$ denote the maximum value of the quality score $Q$ can attain on database $D$. Then, with an $\epsilon$-private exponential mechanism with quality score on database $D$ and outputs satisfies that, with probability at least $1-\beta$, $Q(r, D) \geq \OPT - \frac{2 \Delta}{\epsilon} \log \left(\frac{|\mathcal{R}|}{\beta} \right)$.
\end{lemma}

Finally, we recall a standard composition tool to combine private mechanisms. 

\begin{lemma}[\cite{drv10}] \label{lem:composition_lemma}
Let $\epsilon\in(0,1]$, $\delta \in (0,1/2]$, $\delta'\geq 0$, and let ${\cal M}_1,\ldots,{\cal M}_k$ be $(\epsilon',\delta')$-private, adaptively chosen mechanisms. Then the composition ${\cal M}_1\circ \ldots\circ {\cal M}_k$ is $(\epsilon,\delta+k\delta')$-private, provided that $\epsilon' = \frac{\epsilon}{\sqrt{8k \log (1/\delta)}}$.
\end{lemma}

We remark that while adaptive composition has been improved via moment accountant~\cite{acg+16} and R{\'e}nyi DP~\cite{m17,wbk19}, we use the simpler adaptive composition for illustrating the main idea.
\section{A DP Framework for Euclidean Jordan Algebras}
\label{sec:dp_eja}

In this section, we introduce a differential privacy (DP) framework for Euclidean Jordan algebras (EJA), based on the Gaussian mechanism. We begin with a motivating example. Let ${\cal S}^r$ denote the set of all $r \times r$ real symmetric matrices. Suppose we are designing a private recommendation system, where each client record is represented by an $r$-dimensional feature vector $u$. The goal is to privately release the covariance matrix of the dataset. In this setting, two neighboring databases $D, D'$ differ in exactly one feature vector, and the function $f : {\cal D} \rightarrow {\cal S}^r$ outputs the covariance matrix of the data. It is then natural to define sensitivity using the \emph{matrix Schatten-$p$ norms} rather than entrywise norms. For instance, the $\ell_\infty$ sensitivity corresponds to the spectral norm $\|f(D) - f(D')\|$, while the $\ell_1$ sensitivity corresponds to the nuclear norm $\|f(D) - f(D')\|_{S_1}$. These spectrum-aware notions of sensitivity highlight the inadequacy of mechanisms like the Laplace mechanism, which adds noise entrywise, and motivate the need for a more principled approach. This leads us to consider a general formulation of differential privacy in the context of EJAs, which we refer to as DP-EJA.

\begin{definition}
Let ${\cal D}$ be the universe of databases, and let $f : {\cal D} \rightarrow \J$. For neighboring databases $D, D'$, we define:
\begin{CompactItemize}
\item $f$ has $\ell_1$ sensitivity $\Delta_1$ if $\|f(D) - f(D')\|_1 \leq \Delta_1$;
\item $f$ has $\ell_2$ sensitivity $\Delta_2$ if $\|f(D) - f(D')\|_2 \leq \Delta_2$;
\item $f$ has $\ell_\infty$ sensitivity $\Delta_\infty$ if $\|f(D) - f(D')\|_\infty \leq \Delta_\infty$.
\end{CompactItemize}
\end{definition}

Perhaps the first question one might ask for DP-EJA is: how do we define an additive noise mechanism? It is not immediately clear what constitutes a “Gaussian element” in an EJA $\J$. Thus our first order of business is to develop a Gaussian mechanism for EJAs.

\begin{lemma}[Generic Gaussian mechanism]
\label{lem:gaussian_mech_eja}
Let $\J$ be an EJA with $\dim(\J) = k$, equipped with an isometry $\phi : \J \rightarrow \mathbb{R}^k$. Let $\epsilon, \delta \in (0,1)$ and let $f : \mathcal{D} \rightarrow \J$ be a function with $\ell_2$ sensitivity $\Delta_2$. Consider the following mechanism:
\begin{CompactItemize}
    \item Set $\sigma = \frac{\Delta_2 \sqrt{2 \log(1.25 / \delta)}}{\epsilon}$;
    \item Generate a Gaussian noise vector $\nu \sim \mathcal{N}(0, \sigma^2 I_k)$;
    \item Set $z = \phi^{-1}(\nu)$.
\end{CompactItemize}
The mechanism releases $f(D) + z$ and is $(\epsilon, \delta)$-differentially private.
\end{lemma}

\begin{proof}
The proof is straightforward. Since $\phi$ is an isometry, we have that for any $x, y \in \J$, $\langle \phi(x), \phi(y) \rangle = \langle x, y \rangle$, which implies that $\|x\|_2 = \|\phi(x)\|_2$. Therefore, the $\ell_2$ sensitivity of $f$ satisfies $\|f(D) - f(D')\|_2 = \|\phi(f(D)) - \phi(f(D'))\|_2$. By Definition~\ref{def:gaussian_mech}, the Gaussian mechanism $(\phi \circ f)(D) + \nu$ is $(\epsilon, \delta)$-differentially private. Since $\phi^{-1}$ is also an isometry, we have $\phi^{-1}((\phi \circ f)(D) + \nu) = f(D) + \phi^{-1}(\nu) = f(D) + z$, which implies that the release $f(D) + z$ is $(\epsilon, \delta)$-private.
\end{proof}

Lemma~\ref{lem:gaussian_mech_eja} provides a simple method for constructing a Gaussian element in $\J$: choose an isometry $\phi$, generate a Gaussian vector in $\mathbb{R}^k$, and then apply $\phi^{-1}$. Such an isometry always exists, as one can construct it using an orthonormal basis for $\J$. As an example, consider $\J = {\cal S}^r$, the space of $r \times r$ real symmetric matrices. Then $\dim(\J) = k = \frac{r(r+1)}{2}$, and the isometry $\phi$ maps a symmetric matrix to a vector with the entries from its upper triangular part. Lemma~\ref{lem:gaussian_mech_eja} in this case corresponds to generating a Gaussian vector in $\mathbb{R}^k$, and then applying $\phi^{-1}$ to obtain a symmetric matrix with the upper triangular part filled by the Gaussian vector. The $\ell_2$ norm on ${\cal S}^r$ corresponds to the Frobenius norm, so it is easy to see that Lemma~\ref{lem:gaussian_mech_eja} recovers the standard Gaussian mechanism for matrix-valued functions under Frobenius norm sensitivity. As a consequence, we obtain private mechanisms for $\ell_1$ and $\ell_\infty$ sensitivity via norm inequalities:

\begin{corollary}
\label{cor:l1_eja}
Let $\J$ be an EJA with $\dim(\J) = k$, equipped with an isometry $\phi : \J \rightarrow \mathbb{R}^k$. Let $\epsilon, \delta \in (0,1)$, and let $f : {\cal D} \rightarrow \J$ have $\ell_1$ sensitivity $\Delta_1$. Then, the generic Gaussian mechanism (Lemma~\ref{lem:gaussian_mech_eja}) with sensitivity parameter $\Delta_1$ is $(\epsilon, \delta)$-differentially private.
\end{corollary}

\begin{corollary}
\label{cor:linfty_eja}
Let $\J$ be an EJA with $\dim(\J) = k$, equipped with an isometry $\phi : \J \rightarrow \mathbb{R}^k$. Let $\epsilon, \delta \in (0,1)$, and let $f : {\cal D} \rightarrow \J$ have $\ell_\infty$ sensitivity $\Delta_\infty$. Then, the generic Gaussian mechanism (Lemma~\ref{lem:gaussian_mech_eja}) with sensitivity parameter $\sqrt{r} \Delta_\infty$ is $(\epsilon, \delta)$-differentially private.
\end{corollary}

One might wonder whether it is possible to perturb only the eigenvalues instead of all $k$ dimensions. This idea is particularly tempting for ${\cal S}^r$, since $k = O(r^2)$ while the rank is only $r$. A na\"{i}ve approach would be to first compute a spectral decomposition $f(D) = \sum_{i=1}^r \lambda_i q_i$, then inject scalar Gaussian noise to each eigenvalue, yielding $\sum_{i=1}^r (\lambda_i + \nu_i) q_i$. This approach also appears attractive for handling $\ell_1$ sensitivity by injecting Laplace noise to each eigenvalue. Unfortunately, this method fails to ensure differential privacy. It is not sufficient to perturb only the eigenvalues — the Jordan frame must also be randomized. This can sometimes be achieved by sampling a random Jordan frame $\{p_1, \ldots, p_r\}$ and outputting $\sum_{i=1}^r (\lambda_i + \nu_i) p_i$. However, this raises challenges in bounding the $\ell_2$ norm of the result, as one must account for differences between Jordan frames. In the special case $\J = \mathbb{R}^k$, this issue disappears since there is a unique Jordan frame $\{e_1, \ldots, e_k\}$ and rank coincides with dimension. For $\ell_1$ sensitivity, one might wish to obtain a pure $\epsilon$-private mechanism using the Laplace mechanism: namely, sample a Laplace noise vector in $\mathbb{R}^k$ with suitable parameters and map it into $\J$ via $\phi^{-1}$. Unfortunately, this also fails: the $\ell_1$ norm in $\mathbb{R}^k$ does not, in general, correspond to the $\ell_1$ norm in $\J$ (except when $\J = \mathbb{R}^k$). In essence, the isometry $\phi$ only preserves the inner product — and hence the $\ell_2$ norm — but not other norms. For this reason, our mechanisms for $\ell_1$ and $\ell_\infty$ sensitivity are derived from the generic Gaussian mechanism.

We also note that the utility guarantee for the generic Gaussian mechanism follows from the isometry property of $\phi$: the $\ell_2$ norm of the noise element $z \in \J$ equals the $\ell_2$ norm of the Gaussian vector $\nu \in \mathbb{R}^k$, which scales as $\sigma\sqrt{k}$ due to Lemma~\ref{lem:chi_square_tail}. This implies a weaker bound when translating to $\ell_\infty$ norm. For ${\cal S}^r$, standard results in random matrix theory state that $\|z\|_\infty = O(\sigma\sqrt{r})$~\cite{w58}. For general EJAs, it is well-known that they are direct sums of five simple EJAs~\cite{fk94}, and if it does not have a spin factor as a component, then $\|z\|_\infty=O(\sigma(\sqrt{r_{\max}}+\sqrt{\log k}))$ where $r_{\max}$ is the largest rank among its components, and is $O(\sigma(\sqrt{k_{\max}}+\sqrt{r_{\max}}+\sqrt{\log k}))$ if it has a spin factor as a component and $k_{\max}$ is the largest dimension. To unify the discussion, we conservatively adopt the weaker $O(\sigma\sqrt{k})$ bound in this work.

\section{Private Symmetric Cone Programming}
\label{sec:dp_scp}
In this section, we develop differentially private algorithms for solving symmetric cone programs. As an application, we obtain private algorithms for semidefinite programming, resolving a major open question posed in~\cite{hrru14}. Following the framework of~\cite{hrru14} for linear programming, we study private algorithms under several settings: (1) \emph{High sensitivity constraint privacy}, where neighboring databases may differ by one entire constraint; (2) \emph{Low sensitivity constraint privacy}, where all databases have the same number of constraints, and neighboring databases differ in the $\ell_\infty$ norm of those constraints; (3) \emph{Scalar privacy}, where neighboring databases differ in the right-hand side vector $b$, again measured in $\ell_\infty$ norm; (4) \emph{Objective privacy}, where neighboring instances differ in the objective element $c$ under $\ell_\infty$ norm. Our algorithm for the high sensitivity constraint setting is a generalization of the dense multiplicative weights update (MWU) method used in~\cite{hrru14} for linear programs. For scalar privacy and low sensitivity constraint privacy, we adopt the MWU framework for symmetric cone programming introduced in~\cite{clpv23,zvtl24}, but propose a novel scheme in which the update direction is determined by identifying the most violated constraint. For objective privacy, we perturb the objective element with our generic Gaussian mechanism.

\subsection{High Sensitivity Constraint Privacy}

In this setting, the SCP instances for two neighboring databases share the same objective element $c \in \J$ but differ by one constraint and its corresponding scalar value: they have the same first $m$ constraints, while $D'$ contains an additional constraint and scalar value $b_{m+1}'$. Note that this additional constraint can be any constraint admissible for the program class, for example any $A\succeq 0$ with $\|A\|\leq 1$ for the covering SDPs below. As first observed in~\cite{hrru14}, it is generally impossible to design a private algorithm that satisfies \emph{all} constraints, as a new constraint can significantly alter the optimal solution of the original program. The key idea in~\cite{hrru14} is to run an MWU procedure over a restricted set of constraints, so that the output solution satisfies most of the constraints while preserving privacy. This is achieved by applying Bregman projection~\cite{hw01} and performing MWU over a projected dense distribution over the constraints. While~\cite{hrru14} analyzes this method in the context of linear programming, we show that it extends naturally to symmetric cone programming, since the algorithm operates over constraints. As a consequence, we obtain private algorithms for covering semidefinite programs that are especially useful when the number of constraints $m$ is much larger than $r\cdot\OPT^2/\epsilon$, where $r$ is the matrix dimension and $\OPT$ is the optimal value. This is because our bound on the number of violated constraints grows linearly in $r$ and, up to a logarithmic factor, quadratically in $\OPT/\alpha$ (Theorem~\ref{thm:private_convering_sdp_informal}). Specifically, consider the following covering SDP:
\begin{align*}
    \min_{X \succeq 0} & ~ \Tr(X) \\
    \text{s.t.} & ~ \langle A_i, X \rangle \geq 1, \quad \forall i \in [m],
\end{align*}
where $A_1, \ldots, A_m \succeq 0$ and $\max_{i \in [m]} \|A_i\| \leq 1$. All matrices are of size $r \times r$. Covering SDPs have a wide range of applications in machine learning, including robust mean estimation~\cite{cdg19}, robust covariance estimation~\cite{cdgw19}, and E-optimal experimental design~\cite{vg19}.

\begin{theorem}[Informal version of Theorem~\ref{thm:private_covering_sdp}]
\label{thm:private_convering_sdp_informal}
Let $\epsilon \in (0,1]$, $\delta \in (0,1/2]$ be the DP parameters, let $\beta \in (0,1)$ be the failure probability, and let $\alpha\in(0,\OPT]$ be the target accuracy. Consider a feasible covering SDP with $m\geq 3$ constraints over $r \times r$ matrices, and suppose that its optimal value $\OPT$ and an upper bound $\bar m\geq m+1$ on the number of constraints are given to the algorithm as public inputs, the same on neighboring instances (Remarks~\ref{rem:public_opt} and~\ref{rem:public_rounds}). There exists an algorithm (Algorithm~\ref{alg:constraint_mwu}) that finds $X^* \succeq 0$ with $\Tr(X^*)=\OPT$ such that, with probability at least $1 - \beta$, fewer than $s$ constraints have $\langle A_i, X^* \rangle < 1 - \alpha$, provided that $m\geq s$ and
\begin{align*}
    s \geq C_0\cdot \frac{\OPT^2}{\epsilon\,\alpha^2} \log^{1/2}(1 / \delta) \log^{1/2} \bar m \cdot \Big( r\log\frac{2\OPT}{\alpha} + \log\frac{\log \bar m}{\beta} \Big)
\end{align*}
for an absolute constant $C_0$. Moreover, for every value of the public input $\OPT$ that is at least $1$, the algorithm is $(\epsilon,\delta)$-differentially private with respect to high sensitivity constraint privacy.
\end{theorem}

The core idea behind Theorem~\ref{thm:private_convering_sdp_informal} is to use the dense MWU framework described above, paired with a private oracle that performs a simple linear maximization. To ensure privacy, the oracle is implemented using the exponential mechanism. A major technical challenge is that, even when the optimal value of the program is fixed to $\OPT$, the number of feasible solutions is infinite, as they correspond to extreme rays of the positive semidefinite cone intersected with a hyperplane. A na{\"i}ve application of Lemma~\ref{lem:acc_gua_exp_mec} would yield a vacuous bound due to $|\mathcal{R}| = \infty$. To overcome this, we use a $\gamma$-net argument to discretize the space of feasible solutions. We take $\gamma=\alpha/(12\sqrt{\OPT})$, proportional to the target accuracy. Every point $u$ of the net still gives a point $uu^\top$ of $\K_{\OPT}=\{X\succeq 0:\Tr(X)=\OPT\}$. Rounding to the net changes the oracle's objective by at most $\alpha/6$, and the net has size $\exp(O(r\log(2\OPT/\alpha)))$ (Lemma~\ref{lem:private_oracle}). So Lemma~\ref{lem:acc_gua_exp_mec} applies with a penalty factor of $O(r\log(2\OPT/\alpha))$ in place of $\log|\mathcal{R}|$. This contrasts with the private covering LP algorithm in~\cite{hrru14}, where the candidate solutions are the scaled standard basis vectors in $\mathbb{R}^r$, so $|\mathcal{R}| = r$ and the penalty is only $\log r$. The factor $(\OPT/\alpha)^2$ in $s$, which also appears in the private fractional set cover solver of~\cite{hrru14}, has two sources. First, the score of the exponential mechanism has sensitivity $O(\OPT/s)$ and must be accurate to $O(\alpha)$. Second, multiplicative weights with width $\OPT$ runs for $T=O(\OPT^2\log\bar m/\alpha^2)$ rounds, and by advanced composition the privacy costs of these rounds add up to a factor $\sqrt{T}$. Solving the condition for the target accuracy instead, we find that a given $s$ allows every $\alpha\in(0,\OPT]$ with
\begin{align*}
    \alpha \geq C_1\cdot\OPT\cdot\Big(\frac{r\log(2\OPT/\alpha)+\log(\log\bar m/\beta)}{\epsilon\, s}\Big)^{1/2}\log^{1/4}(1/\delta)\log^{1/4}\bar m
\end{align*}
for an absolute constant $C_1$. Since $\OPT\geq 1$ and $\langle A_i,X\rangle\geq 0$ for every $X\succeq 0$, the guarantee is meaningful only for $\alpha<1$. Our bound reaches this range once $s$ exceeds $r\cdot\OPT^2/\epsilon$ by logarithmic factors.

Inspired by our private algorithm for covering SDP, we further extend the framework to any covering \emph{symmetric cone program} over a simple Euclidean Jordan algebra:
\begin{align*}
    \min_{x \in \K} & ~ \Tr(x) \\
    \text{s.t.} & ~ \langle a_i, x \rangle \geq 1, \quad \forall i \in [m].
\end{align*}

It is well-known that, up to isomorphism, there are five types of simple Jordan algebras: $r \times r$ real symmetric matrices, $r \times r$ complex Hermitian matrices, $r \times r$ quaternionic Hermitian matrices, the spin factors $\R\times\R^n$ with $n\geq 2$, and the exceptional Albert algebra~\cite{fk94}. In particular, the cone of squares of the spin factor $\R\times\R^n$ is the second-order cone $\{(x_0,\bar x)\in\R\times\R^n: x_0\geq\|\bar x\|_2\}$. To apply the machinery developed for private covering SDP to general SCP, we observe that the oracle's optimal solutions are (scaled) primitive idempotents in $\J$. Therefore, a $\gamma$-net argument requires bounding the dimension of the set of primitive idempotents. For all simple Jordan algebras, this dimension is $d(r-1)$, where $d$ is the Peirce constant~\cite{fk94}. The cone $C$ of scaled primitive idempotents then has dimension $k_C=d(r-1)+1$. This dimension is $O(r)$ for the three matrix algebras and $17$ for the Albert algebra. It equals $n$, however, for the spin factor $\R\times\R^n$, whose rank is $2$ (Table~\ref{tab:simple_jordan}). For $\J={\cal S}^r$ it equals $r$, the dimension of the vectors $u$ in the SDP case. To the best of our knowledge, this is the first attempt to quantize the rays of primitive idempotents via a $\gamma$-net, and we hope this approach enables further applications and tighter bounds in future work.

\begin{theorem}[Informal version of Theorem~\ref{thm:private_covering_scp}]
\label{thm:private_convering_scp_informal}
Let $\J$ be a simple EJA of rank $r$ and Peirce constant $d$, and let $k_C=d(r-1)+1$, so that $k_C\leq 4r$ for the three matrix algebras, $k_C=n$ for the spin factor $\R\times\R^n$, and $k_C=17$ for the Albert algebra. Let $\epsilon \in (0,1]$, $\delta \in (0,1/2]$ be the DP parameters, let $\beta \in (0,1)$ be the failure probability, and let $\alpha\in(0,\OPT]$ be the target accuracy. Consider a feasible covering SCP with $m\geq 3$ linear constraints, where $a_i\in\K$ and $\|a_i\|_\infty\leq 1$, and suppose that its optimal value $\OPT$ and an upper bound $\bar m\geq m+1$ on the number of constraints are given to the algorithm as public inputs, the same on neighboring instances (Remarks~\ref{rem:public_opt} and~\ref{rem:public_rounds}). There exists an algorithm (Algorithm~\ref{alg:constraint_mwu}) that outputs a point $x^* \in \K$ with $\Tr(x^*)=\OPT$ such that, with probability at least $1 - \beta$, fewer than $s$ constraints have $\langle a_i, x^* \rangle < 1 - \alpha$, provided that $m\geq s$ and
\begin{align*}
    s \geq C_0\cdot \frac{\OPT^2}{\epsilon\,\alpha^2} \log^{1/2}(1 / \delta) \log^{1/2} \bar m \cdot \Big( k_C\log\frac{2\OPT}{\alpha} + \log\frac{\log \bar m}{\beta} \Big)
\end{align*}
for an absolute constant $C_0$. Moreover, for every value of the public input $\OPT$ that is at least $1$, the algorithm is $(\epsilon,\delta)$-differentially private with respect to high sensitivity constraint privacy.
\end{theorem}

\subsection{Low Sensitivity Constraint Privacy}

In the low sensitivity setting, two neighboring databases have the same number of constraints and differ in the $\ell_\infty$ norm. Specifically, let $A(D), A(D')$ be the constraint sets corresponding to neighboring databases $D$ and $D'$, respectively. We assume that $\max_{i \in [m]} \| a_i(D) - a_i(D') \|_\infty \leq \Delta_\infty$. As observed in~\cite{hrru14} for the LP setting, it is possible to approximately satisfy all constraints by applying MWU over the \emph{variables}. In particular, the oracle in this setting is given the constraint set $A \in \J^m$, scalar vector $b \in \mathbb{R}^m$, and a point $x \in \K$, and it must return an approximately most violated constraint $a_i$ such that $\langle a_i, x \rangle - b_i \geq \max_{j \in [m]} \langle a_j, x \rangle - b_j - \alpha$ with high probability. We refer to such an oracle as an \emph{$(\alpha, \gamma)$-dual oracle}, which achieves additive approximation $\alpha$ and fails with probability at most $\gamma$. This oracle is particularly well-suited to privatization via the exponential mechanism: the score function $Q(i, b) = \langle a_i, x \rangle - b_i$ leads to both privacy and accuracy guarantees. It remains to show that one can indeed solve the SCP in a first-order fashion using such an oracle. In the case of linear programming, this is already established by the classical work of~\cite{pst95}. Our first result shows that an approximate most violated constraint oracle can likewise be used to solve SCPs over symmetric cones in a first-order manner. Before stating the result, we define the \emph{width} of the constraint set as $\rho = \max_{i \in [m]} \| a_i \|_\infty$.

\begin{theorem}[Informal version of Theorem~\ref{thm:non_private_scp}]
\label{thm:non_private_scp_informal}
Given an SCP with $m$ linear constraints, suppose there exists a feasible point $x \in \K$ with $\Tr(x) = 1$, and access to an $(\alpha/2, \gamma)$-dual oracle. Then, there exists an algorithm (Algorithm~\ref{alg:primal_lp_nonpriv}) that finds a distribution element $x^* \in \K$ such that $\langle a_i, x^* \rangle \leq b_i + \alpha$ for all $i \in [m]$, with probability at least $1 - T \gamma$, where $T = O\left( \frac{\rho^2 \log r}{\alpha^2} \right)$.
\end{theorem}

For a program with a feasible solution of trace $L>0$, the assumption on the existence of $x$ with $\Tr(x) = 1$ can be met by dividing that solution and the vector $b$ by $L$ (Appendix~\ref{sec:low_sensitivity}). Our privacy guarantees then need $L$ to be public, the same on neighboring instances, since a scale computed from the private data would not be private. Previous MWU algorithms for symmetric cone programming have relied on a \emph{primal oracle}, where the input is a point $x \in \K$ with $\Tr(x) = 1$, and the output is a vector $y \in \mathbb{R}^m$ satisfying $\langle \sum_i y_i a_i - c, x \rangle \geq 0$ and $b^\top y \leq \alpha$~\cite{zvtl24}. Our proof instead builds on the regret bound for MWU over EJAs as established in~\cite{clpv23}.

\begin{algorithm}[ht!]\caption{Constraint private SCP solver.}
  \label{alg:primal_lp_rowsens}
  \begin{algorithmic}[1]
    \Procedure{ConstraintPrivateSCP}{$A\in \J^m, b\in \R^m$}
    \State $x^1\gets e/r$
    \State Let $\alpha,\gamma$ be the parameters for the oracle, $\epsilon,\delta$ be the parameters for DP
    \State Let $\textsc{Oracle}$ be an $(\alpha,\gamma)$-dual oracle
    \State $T\gets \frac{144\log r}{\alpha^2}, \epsilon'\gets \frac{\epsilon}{4\sqrt{T\log(2/\delta)}}, \eta\gets \frac{\alpha}{12}$
    \For{$t=1\to T$}
    \State $\sigma\gets \frac{\Delta_\infty\sqrt{2r\log(5T/\delta)}}{\epsilon'}$
    \State $p^t\gets \textsc{Oracle}(A, b, x^t)$
    \State $\nu^t\sim \N(0, \sigma^2 I_k)$
    \State $z^t\gets \phi^{-1}(\nu^t)$
    \State $\wh \ell^t\gets \frac{a_{p^t}+z^t}{2}$
    \State  $x^{t+1}\gets \frac{\exp(-\sum_{i=1}^t \eta \wh \ell^i)}{\Tr(\exp(-\sum_{i=1}^t \eta \wh \ell^i))}$
    \EndFor
    \State \Return $\ov x\gets \frac{1}{T}\sum_{t=1}^T x^t$
    \EndProcedure
  \end{algorithmic}
\end{algorithm}
We observe that the private data $A$ is accessed both during the oracle step and the loss computation. To ensure privacy in the oracle step, we use the exponential mechanism; for the loss computation, we apply the generic Gaussian mechanism developed in Section~\ref{sec:dp_eja}, injecting a Gaussian noise into the constraint returned by the oracle.

\begin{theorem}[Informal version of Theorem~\ref{thm:rowsens_priv} and~\ref{thm:rowsens_acc}]
\label{thm:rowsens_informal}
Let $A \in \J^m$ satisfy $\lambda(a_i) \subseteq [-1, 1]$ for all $i \in [m]$, and let $b \in \mathbb{R}^m$. Suppose the SCP admits a feasible $y \in \K$ with $\Tr(y) = 1$. Let $\beta \in (0,1)$, $\epsilon \in (0,1]$ and $\delta \in (0, 1/2]$. Then Algorithm~\ref{alg:primal_lp_rowsens} using the exponential mechanism to implement a dual oracle returns a distributional element $x^*$ such that with probability at least $1-\beta$, $\langle a_i,x^*\rangle\leq b_i+\alpha$ for all $i \in [m]$, where
\begin{align*}
    \alpha = \widetilde{O}\left( \frac{\Delta_\infty^{1/2} r^{1/4} k^{1/4}}{\epsilon^{1/2}} \cdot \mathrm{polylog}(r, m, 1/\beta, 1/\delta) \right).
\end{align*}
Moreover, the algorithm is $(\epsilon, \delta)$-private with respect to low sensitivity constraint privacy.
\end{theorem}

Compared to the low sensitivity constraint privacy LP result of~\cite{hrru14}, our bound incurs an additional $k^{1/4}$ factor in $\alpha$. This arises because the sensitivity is measured in the $\ell_\infty$ norm, yet the noise added is a $k$-dimensional Gaussian vector. By standard concentration bounds on the norm of Gaussian vectors (see Lemma~\ref{lem:chi_square_tail}), the $\ell_2$ norm of the Gaussian vector — and consequently, the corresponding Gaussian element $z$ — scales with $\sqrt{k}$. In contrast,~\cite{hrru14} adds independent entrywise Laplace noise, and the magnitude of any single perturbation is $O(\Delta_\infty / \epsilon)$ up to a logarithmic factor with high probability. As discussed in Section~\ref{sec:dp_eja}, our framework does not permit injecting noise directly into the eigenvalues of the EJA elements, since this would require sampling a random Jordan frame as well — posing additional complexity and potential distortion. We also develop algorithms for the setting where the scalar vector $b \in \mathbb{R}^m$ or the objective element $c \in \J$ is private under low sensitivity. For scalar privacy, we run Algorithm~\ref{alg:primal_lp_nonpriv} with the exponential mechanism as the dual oracle. For objective privacy, we perturb $c$ with the generic Gaussian mechanism and solve the perturbed program with the added constraint $\|x\|_2=1$, assuming an optimal solution of unit $\ell_2$ norm. We defer the details to Appendix~\ref{sec:low_sensitivity}.
\begin{remark}
We interpret Algorithm~\ref{alg:primal_lp_rowsens} as an approximate multiplicative weights update (MWU) scheme with noisy oracles. Instead of applying the MWU rule directly to the primal variables, the algorithm first perturbs the oracle output with Gaussian noise and then uses the perturbed value for the update. Theorem~\ref{thm:rowsens_informal} can thus be viewed as quantifying how the injected noise influences constraint violations. More broadly, by altering the noise distribution, one can adapt the noisy MWU framework to other settings that require different types of error guarantees.
\end{remark}
\section{Conclusion}
\label{sec:concl}

We study differentially private algorithms for convex programming, with a particular focus on symmetric cone programming. To this end, we develop a generic Gaussian mechanism for Euclidean Jordan algebras that provides differential privacy guarantees under $\ell_1$, $\ell_2$, and $\ell_\infty$ norms. We incorporate this mechanism into our private solvers for low sensitivity constraint privacy and objective privacy. Our solvers for scalar privacy and for the high sensitivity setting use the exponential mechanism instead. For the high sensitivity regime, we generalize the analysis of~\cite{hrru14} beyond linear programming and apply it to covering semidefinite programs and covering symmetric cone programs. In the low sensitivity setting, we design private solvers for scalar and constraint privacy based on an approximately most violated constraint oracle. As a direct consequence, we obtain a family of differentially private algorithms for semidefinite programming — a longstanding open problem originally posed by~\cite{hrru14} — that also encompass a wide range of applications in machine learning.

There are several limitations of our work, which we leave as directions for future research. (1) The differential privacy mechanisms we employ are relatively basic; recent advances in moment accounting and R{\'e}nyi DP~\cite{acg+16,m17,wbk19} could yield stronger trade-offs between privacy and utility. (2) Although our private solvers provide meaningful guarantees, all of them except the objective private one only approximately satisfy the constraints. In high sensitivity settings, this manifests as a small number of constraint violations. For low sensitivity scalar and constraint privacy, a central open question is whether it is possible to design private algorithms that satisfy all constraints exactly. For linear programming,~\cite{muvv21,bbdh24,bbdh25} develop techniques to preserve feasibility under privacy. Extending these methods to the broader setting of symmetric cone programming (and even SDPs) is considerably more challenging, since one must contend with privacy-preserving perturbations to the \emph{spectrum} of matrices rather than to entries of a vector. (3) In the high sensitivity setting, our bound $s$ on the number of violated constraints grows as $(\OPT/\alpha)^2$ up to a logarithmic factor, as in the private fractional set cover solver of~\cite{hrru14}. Later work~\cite{cz26} improved the dependence on the scale and the accuracy from polynomial to logarithmic, at the price of a bound that grows as $k^2$ in the dimension $k$ of the ambient space, where our bound for covering SDPs grows linearly in the matrix dimension $r$ (see the paragraph on subsequent work in Section~\ref{sec:relat}). Our covering solvers also take $\OPT$, or a public upper bound on it, as input. We leave open whether a bound can be linear in $r$ and logarithmic in $\OPT/\alpha$ at the same time, and whether a private estimate of $\OPT$ can replace the public value at a small cost.

\ifdefined\isarxiv
\section*{Acknowledgment}

The corrections in this version were found and written with the help of Claude Opus 5.5. The authors have verified the corrections and take full responsibility for the content of this paper.
\else
\section*{Acknowledgment}

We thank anonymous NeurIPS reviewers for their constructive comments. Lichen Zhang is supported by a Mathworks Fellowship and a Simons Dissertation Fellowship in Mathematics.
\fi

\ifdefined\isarxiv
\bibliographystyle{alpha}
\bibliography{ref}
\else
\bibliographystyle{alpha}
\bibliography{ref}
\fi


\newpage
\onecolumn
\ifdefined\isarxiv
\else
\clearpage
\clearpage
\section*{NeurIPS Paper Checklist}

\begin{enumerate}

\item {\bf Claims}
    \item[] Question: Do the main claims made in the abstract and introduction accurately reflect the paper's contributions and scope?
    \item[] Answer: \answerYes{} 
    \item[] Justification: The main contributions of this paper include a generic Gaussian mechanism for Euclidean Jordan algebra and private algorithms for symmetric cone programming. They are clearly reflected in the abstract and introduction.

\item {\bf Limitations}
    \item[] Question: Does the paper discuss the limitations of the work performed by the authors?
    \item[] Answer: \answerYes{} 
    \item[] Justification: The paper clearly discusses the limitations in the second paragraph of Section~\ref{sec:concl}.

\item {\bf Theory assumptions and proofs}
    \item[] Question: For each theoretical result, does the paper provide the full set of assumptions and a complete (and correct) proof?
    \item[] Answer: \answerYes{} 
    \item[] Justification: All assumptions and proofs are provided, specifically, see Section~\ref{sec:dp_eja} for the privacy proof of the generic Gaussian mechanism, Appendix~\ref{sec:high_sensitiviy} for the settings and proofs for high sensitivity constraint privacy with applications, and Appendix~\ref{sec:low_sensitivity} for a slew of results in the low sensitivity setting.

    \item {\bf Experimental result reproducibility}
    \item[] Question: Does the paper fully disclose all the information needed to reproduce the main experimental results of the paper to the extent that it affects the main claims and/or conclusions of the paper (regardless of whether the code and data are provided or not)?
    \item[] Answer: \answerNA{} 
    \item[] Justification: The paper is theoretical in nature and does not contain experiments.

\item {\bf Open access to data and code}
    \item[] Question: Does the paper provide open access to the data and code, with sufficient instructions to faithfully reproduce the main experimental results, as described in supplemental material?
    \item[] Answer: \answerNA{} 
    \item[] Justification: The paper is theoretical and does not contain data and code.

\item {\bf Experimental setting/details}
    \item[] Question: Does the paper specify all the training and test details (e.g., data splits, hyperparameters, how they were chosen, type of optimizer, etc.) necessary to understand the results?
    \item[] Answer: \answerNA{} 
    \item[] Justification: The paper is theoretical hence contains no experiments.

\item {\bf Experiment statistical significance}
    \item[] Question: Does the paper report error bars suitably and correctly defined or other appropriate information about the statistical significance of the experiments?
    \item[] Answer: \answerNA{} 
    \item[] Justification: The paper is theoretical hence contains no experiments.

\item {\bf Experiments compute resources}
    \item[] Question: For each experiment, does the paper provide sufficient information on the computer resources (type of compute workers, memory, time of execution) needed to reproduce the experiments?
    \item[] Answer: \answerNA{} 
    \item[] Justification: The paper is theoretical hence contains no experiments.
    
\item {\bf Code of ethics}
    \item[] Question: Does the research conducted in the paper conform, in every respect, with the NeurIPS Code of Ethics \url{https://neurips.cc/public/EthicsGuidelines}?
    \item[] Answer: \answerYes{} 
    \item[] Justification: Authors have reviewed and ensured the paper adheres to the NeurIPS Code of Ethics.

\item {\bf Broader impacts}
    \item[] Question: Does the paper discuss both potential positive societal impacts and negative societal impacts of the work performed?
    \item[] Answer: \answerYes{} 
    \item[] Justification: Broader societal impacts have been discussed in Appendix~\ref{sec:impact}. 
    
\item {\bf Safeguards}
    \item[] Question: Does the paper describe safeguards that have been put in place for responsible release of data or models that have a high risk for misuse (e.g., pretrained language models, image generators, or scraped datasets)?
    \item[] Answer: \answerNA{} 
    \item[] Justification: The paper is theoretical in nature, therefore does not contain data or models.

\item {\bf Licenses for existing assets}
    \item[] Question: Are the creators or original owners of assets (e.g., code, data, models), used in the paper, properly credited and are the license and terms of use explicitly mentioned and properly respected?
    \item[] Answer: \answerNA{} 
    \item[] Justification: The paper is theoretical in nature, therefore does not use any assets.

\item {\bf New assets}
    \item[] Question: Are new assets introduced in the paper well documented and is the documentation provided alongside the assets?
    \item[] Answer: \answerNA{} 
    \item[] Justification: The paper is theoretical in nature, therefore does not create any new assets.

\item {\bf Crowdsourcing and research with human subjects}
    \item[] Question: For crowdsourcing experiments and research with human subjects, does the paper include the full text of instructions given to participants and screenshots, if applicable, as well as details about compensation (if any)? 
    \item[] Answer: \answerNA{} 
    \item[] Justification: The paper is theoretical in nature, therefore does not have any crowdsourcing or research with human subjects.

\item {\bf Institutional review board (IRB) approvals or equivalent for research with human subjects}
    \item[] Question: Does the paper describe potential risks incurred by study participants, whether such risks were disclosed to the subjects, and whether Institutional Review Board (IRB) approvals (or an equivalent approval/review based on the requirements of your country or institution) were obtained?
    \item[] Answer: \answerNA{} 
    \item[] Justification: The paper is theoretical in nature, therefore does not require and IRB approvals.

\item {\bf Declaration of LLM usage}
    \item[] Question: Does the paper describe the usage of LLMs if it is an important, original, or non-standard component of the core methods in this research? Note that if the LLM is used only for writing, editing, or formatting purposes and does not impact the core methodology, scientific rigorousness, or originality of the research, declaration is not required.
    \item[] Answer: \answerNA{} 
    \item[] Justification: LLM was only used for word editing.

\end{enumerate}
\fi
\newpage
\appendix

\begin{center}
    \textbf{\LARGE Appendix }
\end{center}

\section{More Preliminaries}
\label{sec:more_preli}
In this section, we provide more preliminaries on standard inequalities for EJA and multiplicative weights update for SCP.

Cauchy-Schwarz inequality is automatically satisfies given an proper inner product, we record it here.
\begin{definition}[Cauchy-Schwarz inequality for EJA]
Let $x, y\in \J$, then we have the following standard Cauchy-Schwarz inequality:
\begin{align*}
    |\langle x,y\rangle|\leq & ~ \|x\|_2\cdot \|y\|_2.
\end{align*}
\end{definition}

While Cauchy-Schwarz inequality only requires an inner product space, H{\"o}lder's inequality requires dual norms and a pairing. Nevertheless, we prove the H{\"o}lder's inequality for EJA.

\begin{lemma}[H{\"o}lder's inequality for EJA]\label{lem:holder}
Let $p,q\in [1,\infty]$ satisfy $1/p+1/q=1$, let $x,y\in \J$, then we have
\begin{align*}
    |\langle x, y\rangle| \leq & ~ \|x\|_p\cdot \|y\|_q.
\end{align*}
\end{lemma}

\begin{proof}
The idea is to write $x, y$ in their spectral decomposition and apply the standard vector H{\"o}lder's inequality. Let $x=\sum_{i=1}^r \lambda_i u_i, y=\sum_{i=1}^r \mu_i v_i$ be their respective spectral decomposition. Note that given a primitive idempotent element $u$, we have that $\Tr(u)=\|u\|_2^2=1$ because its rank-1. Then
\begin{align*}
    |\langle x,y\rangle| = & ~ |\Tr(x\circ y)| \\
    = & ~ |\Tr(\sum_{i=1}\lambda_i u_i\circ \sum_i \mu_i v_i)| \\
    = & ~ |\sum_{i,j} \lambda_i\mu_j \Tr(u_i\circ v_j)| \\
    \leq & ~ \sum_{i=1}^r \sum_{j=1}^r |\lambda_i|\cdot |\mu_j| \cdot \langle u_i,v_j\rangle,
\end{align*}
consider the matrix $P_{i,j}=\langle u_i,v_j\rangle$, note that $P$ is a doubly stochastic matrix: since all $u_i$'s, $v_j$'s are primitive idempotent, it must be that $\langle u_i,v_j\rangle\geq 0$. Moreover, it's easy to verify that if we fix $i$, then
\begin{align*}
    \sum_{j=1}^r \langle u_i,v_j\rangle = & ~ \langle u_i,\sum_{j=1}^r v_j\rangle \\
    = & ~ \langle u_i,e\rangle \\
    = & ~ \Tr(u_i) \\
    = & ~ 1,
\end{align*}
and the same argument holds for rows. By Birkhoff–von Neumann theorem~\cite{b46}, $P$ can be written as a convex combination of permutation matrices :$P=\sum_{i=1}^r c_i\Pi_i$ where $c_i$'s form a convex combination, so it suffices to work with any permutation matrix $\Pi$ and its corresponding permutation function $\sigma:[r]\rightarrow [r]$, set $|\lambda|,|\mu|$ be the vector with entries in $|\lambda_i|, |\mu_j|$, then
\begin{align}\label{eq:sorting}
    \sum_{i,j} |\lambda_i|\cdot |\mu_j|\cdot \Pi_{i,j} = & ~ |\lambda|^\top \Pi |\mu| \notag\\
    = & ~ \sum_{i=1}^r |\lambda_i|\cdot |\mu_{\sigma(i)}|.
\end{align}
Thus, we can conclude that
\begin{align*}
    \sum_{i,j} |\lambda_i|\cdot |\mu_j|\cdot \langle u_i,v_j\rangle = & ~ |\lambda|^\top P |\mu| \\
    = & ~ \sum_{i=1}^r c_i |\lambda|^\top \Pi_i |\mu| \\
    \leq & ~ \sum_{i=1}^r c_i \| \lambda\|_p \|\Pi_i \mu\|_q \\
    = & ~ \sum_{i=1}^r c_i \|\lambda\|_p\|\mu\|_q \\
    = & ~ \|\lambda\|_p\|\mu\|_q \\
    = & ~ \|x\|_p\|y\|_q,
\end{align*}
where we use vector H{\"o}lder's inequality in the third step. This completes the proof.
\end{proof}

\cite{twk22} proves a generalized Golden-Thompson inequality for EJA, as follows.
\begin{lemma}[Generalized Golden-Thompson inequality, \cite{twk22}]\label{lem:golden_thompson_inequality}
Let $(\J, \circ)$ be an EJA and $x,y \in \J$, the generalized Golden-Thompson inequality holds:
\begin{align*}
    \Tr(\exp(x+y)) \leq \Tr(\exp(x) \circ \exp(y)).
\end{align*}
\end{lemma}

We now describe the online linear optimization framework over symmetric cones.

\begin{definition}[Online linear optimization (OLO) framework using symmetric cone multiplicative weights update (SCMWU),~\cite{clpv23}]\label{def:online_linear_optimization_scmwu}
    At each time $t$, the algorithm needs to pick a distributional element $p^{t}$ such that $\Tr(p^{t})=1$ and $p^{t}\in \K$. Once the element is picked, a linear loss function $\ell^{t}(p)=\langle m^{t}, p\rangle$ where $m^{t}\in \J$. The SCMWU computes the next iterate as
    \begin{align*}
        p^{t+1} = \frac{\exp(-\eta\sum_{\tau = 1}^t m^{\tau})}{\Tr(\exp(-\eta\sum_{\tau = 1}^t m^{\tau}))},
    \end{align*}
    where $p^{1}=e/\Tr(e)$ is the uniform distribution over $\K$.
\end{definition}

The main idea of SCMWU is to incur cumulative losses comparable to the losses of the best set of actions that an algorithm can make in hindsight. The discrepancy between the cumulative losses of our algorithm and the optimal algorithm is usually referred to as the regret as a function of time $t$. We recall the regret bound of SCMWU proved in~\cite{clpv23}.

\begin{theorem}[Theorem 5.1 of~\cite{clpv23}]
\label{thm:mw_regret}
Let $({\cal J},\circ)$ be an EJA of rank $r$, ${\cal K}$ be its cone of squares. For any $\eta\in (0,1]$ and any sequence of loss vectors $\{\ell^1,\ldots,\ell^T\}$ satisfy $\|\ell^t\|_\infty\leq 1$, the iterates $A^t$ generated by Algorithm~\ref{alg:mw} satisfy
\begin{align*}
    \sum_{t=1}^T \langle \ell^t, A^t\rangle \leq & ~ \sum_{t=1}^T\langle \ell^t, B \rangle+\eta T+\frac{\ln r}{\eta},
\end{align*}
where $B$ is any point in ${\cal K}$ satisfying $\Tr(B)=1$.
\end{theorem}

\section{High Sensitivity Constraint Private SCP}
\label{sec:high_sensitiviy}

In this section, we study algorithms for SCP in the high sensitivity constraint privacy setting. We generalize the algorithm and analysis of~\cite{hrru14} in the LP setting. 
\subsection{Solving SCP with Dense Multiplicative Weights Update}\label{sec:solve_scp_dense_scmwu}
Let us begin by considering constraint private SCPs over a symmmetric cone $\K$ in an EJA $\J$, with the general form
\begin{align*}
  \max_{x \in \mathcal{J}} & ~ \langle c, x \rangle \\
  \text{s.t.} & ~ \langle a_j, x \rangle \leq b_j, \forall j \in [m] \\
  & ~ x \in \mathcal{K},
\end{align*}
where $c, a_1, \ldots, a_m \in \J, b_1, \ldots, b_m \in \R$ and ${\cal K}\subseteq {\cal J}$.

Let $\K_{\OPT} = \mathcal{K} \cap \{ x \in \J \mid \langle c, x \rangle = \OPT \}$. We reduce the SCP to the feasibility program
\begin{align*}
  \text{find } & ~ x  \in \mathcal{K}_{\OPT} \\
  \text{s.t. } & ~ \langle a_j, x \rangle \leq b_j, \forall j \in [m],
\end{align*}
where the value $\OPT$ is a public input to the algorithm, the same on neighboring instances. The privacy guarantees below hold for every fixed value of this input. A binary search for $\OPT$ would read the private constraints, and our privacy analysis does not cover it. Hence, it suffices to solve the feasibility program. As $\K_{\OPT}$ is convex, so we write $\K$ for $\K_{\OPT}$ for simplicity.

Let $A = \{a_1, \ldots, a_m\} \in \J^m$ be a collection of constraints. A database $D$ defines an SCP as a tuple $(c(D), A(D), b(D))$, where $c(D)$ is the objective element, $A(D)$ is the set of constraints and $b(D)$ is the right hand vector. In a constraint private SCP, we have that $c(D)=c(D')$, the only differing parts are the constraints $A(D), A(D')$ and their associated scalar value. In particular, $|A(D)\cap  A(D')|=m, |A(D)|=m, |A(D')|=m+1$, so $A(D), A(D')$ coincide with all $m$ but one constraint. The scalars for the coinciding $m$ constraints are the same, except the scalar value for the differing constraint. We define the adjacency as differing by exactly one constraint.

\begin{definition}[High sensitivity constraint privacy]\label{def:high_sensitivity_constraint_privacy}
Given $m \in \mathbb{N}$, vector $b \in \mathbb{R}^m$, and a constraint set $A \in \J^m$, 
a randomized mechanism $\mathcal{M}$ that outputs a vector in $\J$ is $(\epsilon, \delta)$-high sensitivity constraint private if for any $A, A'$ such that $A'=A\cup \{a_{m+1}\}$, and $b, b'$ such that $b'=[b; b_{m+1}]^\top$,
\begin{align*}
    \Pr[\mathcal{M}(m, b, A) \in S] \leq & ~ e^\epsilon \Pr[\mathcal{M}(m+1, b', A') \in S] + \delta \quad\text{and} \\
    \Pr[\mathcal{M}(m+1, b', A') \in S] \leq & ~ e^\epsilon \Pr[\mathcal{M}(m, b, A) \in S] + \delta
\end{align*}
for any subset $S \subseteq \J$.
\end{definition}

Next, we introduce the dense MWU framework. Let ${\cal F}$ denote the universe of actions, $F$ be the measures on the set of actions $F:{\cal F}\rightarrow [0,1]$, and $\wt F$ be the respective probability distribution, defined as $\wt F=\frac{F}{|F|}$ where $|F|=\sum_{f\in {\cal F}} F_f$ is density of $F$. A key concept we will be relying on is the \emph{Bregman projection}.
\begin{definition}[Bregman projection]\label{def:bregman_projection}
Let $F$ be a measure with $|F|\leq s$ for some $s>0$, we define $\Gamma_s F$ as the Bregman projection of $F$ onto the set of $1/s$-dense distributions:
\begin{align*}
    \Gamma_s F_f := \frac{1}{s} \cdot \min\{1, c F_f\} \quad \forall f \in {\cal F},
\end{align*}
where $c \geq 0$ is the value satisfying $s = \sum_{f \in {\cal F}} \min\{1, c F_f\}$.
\end{definition}

\begin{algorithm}[!ht]
\caption{Dense multiplicative weights update,~\cite{hw01}.}
\label{alg:dense_mwu}
\begin{algorithmic}[1]
\Procedure{DenseMWU}{${\cal F}, \eta$}
\State $F^1 \gets \text{the uniform distribution on ${\cal F}$}$ \Comment{$F^1\in \R^{|{\cal F}|}$}
\For{$t=1\to T-1$}
\State $B^t\gets \Gamma_s F^t$
\State Receive loss vector $\ell^t$ \Comment{$\ell^t\in \R^{|{\cal F}|}$}
\For{$f\in {\cal F}$}
\State $F^{t+1}_f\gets e^{-\eta \ell^t_f} F^t_f$
\EndFor
\EndFor
\State \Return $F^T$
\EndProcedure
\end{algorithmic}
\end{algorithm}

The following lemma due to~\cite{hw01} gives the regret bound for Algorithm~\ref{alg:dense_mwu}.

\begin{lemma}[\cite{hw01}]\label{lem:scp_mw_bound}
Let $F_1$ denote the uniform distribution over $\mathcal{F}$ (so that $|F_1|=1$). 
Consider the sequence of projected distributions $\{\widetilde{B}^t\}_{t=1}^T$ produced by 
Algorithm~\ref{alg:dense_mwu} under an arbitrary loss sequence $\{\ell^t\}_{t=1}^T$ with 
$\|\ell^t\|_\infty \le 1$ and step size $\eta \le 1/2$. 
Define $\widetilde{B}^*$ as the uniform distribution supported on some subset 
$S^* \subseteq \mathcal{F}$ of size $s$.
 Then,
\begin{align*}
    \frac{1}{T} \sum_{t=1}^T \langle \ell^t, \wt{B}^t \rangle \leq \frac{1}{T} \sum_{t=1}^T \langle \ell^t, \wt{B}^* \rangle + \eta + \frac{\log |\mathcal{F}|}{\eta T}.
\end{align*}
\end{lemma}

As standard for any MWU algorithm, we define the width of an SCP, denoted as $\rho$:
\begin{align*}
    \rho \geq \max_D \max_{x \in \mathcal{K}} \max_{i \in [m] } |\langle a_i(D), x \rangle - b_i(D)|,
\end{align*}

MWU algorithm usually assumes an oracle that can efficiently perform ``simpler'' minimization problem, and in our case, it corresponds to minimize over a simple linear function.

\begin{definition}[$(\alpha,\beta)$-approximate, $\rho$-bounded oracle]
Given a distribution $y \in \mathbb{R}^m$ and a set $A=\{a_1,\ldots,a_m\}\in\mathcal{J}^m$, an $(\alpha,\beta)$\textit{-approximate, $\rho$-bounded oracle} returns $x^*\in\mathcal{K}$ with probability at least $1-\beta$ such that
\begin{align*}
\langle \sum_{i=1}^m y_i a_i,\, x^*\rangle \leq \min_{x\in\mathcal{K}} \langle \sum_{i=1}^m y_i a_i,\, x\rangle + \alpha
\quad\text{and}\quad
\max_{i\in[m]} \big|\langle a_i, x^*\rangle - b_i\big| \le \rho.
\end{align*}
\end{definition}

To solve symmetric cone programs, we employ the dense symmetric cone multiplicative weights update algorithm, which maintains a distribution over the constraints and, at each iteration, selects a point $x^t \in \mathcal{K}$ that approximately minimizes the weighted violation of those constraints. Intuitively, losses increase the weight assigned to violated constraints, thereby steering subsequent iterates toward improved feasibility. Averaging the points $x^t$ across all iterations then produces an approximately feasible solution. The full procedure is described in Algorithm~\ref{alg:constraint_mwu}.

\begin{algorithm}[!ht]
\caption{SCP feasibility via \textsc{DenseSCMWU}.}
\label{alg:constraint_mwu}
\begin{algorithmic}[1]
\Procedure{SCPDenseMWU}{$A, b$}
\State $\wt y^1\gets {\bf 1}_m/m$
\State Let $
      \rho \geq \max_D \max_{x \in \mathcal{K}} \max_{i \in [m] } |\langle a_i(D), x \rangle - b_i(D)|$ be the width of the SCP, $s \in \mathbb{N}$ be the density parameter, $\alpha > 0$ be the
      desired accuracy
\State{Let \textsc{Oracle} be an $(\alpha, \beta)$-accurate, $\rho$-bounded oracle}

\State{$T \gets 36 \alpha^{-2} \rho^2 \log m$, $\eta \gets \sqrt{(\log m)/T}$, with $m$ read as a public bound $\bar m\geq m+1$ on the number of constraints (Remark~\ref{rem:public_rounds})}

\For{$ t=1\to T$}
\State{$x^{t} \gets \textsc{Oracle}(\wt{y}^t, A)$}
\State{$\ell^t_i \gets (1/2\rho) (b_i - \langle a_i, x^t\rangle) +
      1/2$}
\State{Update $\wt{y}^{t+1}$ from $\wt{y}^t$ and $\ell^t$ via
      dense multiplicative weights with density $s$}
\EndFor
\State \Return $\ov{x} \gets (1/T) \sum_{t = 1}^T x^t$
\EndProcedure
\end{algorithmic}
\end{algorithm}

Algorithm~\ref{alg:constraint_mwu} differs significantly from other MWU algorithms for solving SCP such as~\cite{clpv23,zvtl24}, where the algorithm updates the \emph{variable $x\in \K$} instead of a distribution over the constraints. In the following, we show that Algorithm~\ref{alg:constraint_mwu} provides utility guarantee for solving the program. The proof is similar to~\cite{hrru14}, as the analysis focuses on the constraint, which does not exploit the structure of EJA $\J$.

\begin{lemma}[Approximate SCP feasibility via \textsc{SCPDenseMWU}]\label{lem:approximate_feasibility}
Let $A=\{a_1,\ldots,a_m\}\in \J^m$ and $b\in \R^m$, $\rho$ be the width of the SCP. Let $\alpha\in (0,3\rho]$, $\beta\in (0,1)$ and $T=36\alpha^{-2}\rho^2\log m$. Suppose the SCP is feasible, then Algorithm~\ref{alg:constraint_mwu} with density $s$ that utilizes an $(\alpha/3,\beta/T)$-approximate, $\rho$-bounded oracle can output $x^*$ in ${\cal K}$ such that, with probability at least $1-\beta$, there exists a subset of constraints $S\subseteq [m]$ with $|S|<s$ such that $\langle a_i,x^*\rangle\leq b_i+\alpha$ for all $i\in [m]\setminus S$.\end{lemma}
\begin{proof}
We will condition on the event that the oracle succeeds on all steps, note that this could be achieved via a union bound over $T$ steps to obtain a success probability of at least $1-\beta$.

    Let $\mathcal{K}_s = \{y \in \R^m \mid \langle {\bf 1}_m, y\rangle =1, \|y\|_\infty \leq 1/s\}$ be the set of $1/s$-dense distribution. The oracle finds $x^t$ with $\sum_{i=1}^m y_i\langle a_i,x^{t}\rangle\leq \sum_{i=1}^m y_i b_i + \alpha/3$. Define the $i$-th item in the loss vector is $\ell^t_i = (1/2\rho)(b_i - \langle a_i,x^t_i\rangle)+1/2$. Then we have
    \begin{align*}
        \langle \ell^t, y^t\rangle = & ~ \sum_{i=1}^m \ell^t_i y_i^t \\
        = & ~ (1/2\rho)\sum_{i=1}^m (b_i-\langle a_i,x_i^t\rangle)y_i^t+\frac{1}{2} \|y\|_1 \\
        \geq & ~ (1/2\rho)((\sum_{i=1}^m y_i^t b_i-y_i^tb_i)-\alpha/3) + \frac{1}{2}\|y\|_1 \\
        = & ~ \frac{1}{2}-\frac{\alpha}{6\rho},
    \end{align*}
    which is nonnegative as $\alpha\leq 3\rho$. To prove an upper bound, we recall that the oracle is $\rho$-bounded:
    \begin{align*}
        \langle \ell^t, y^t\rangle = & (1/2\rho)\sum_{i=1}^m (b_i-\langle a_i,x_i^t\rangle)y_i^t+\frac{1}{2} \|y\|_1 \\
        \leq & ~ \frac{1}{2\rho}\sum_{i=1}^m |b_i-\langle a_i,x_i^t\rangle| y_i^t+\frac{1}{2} \\
        \leq & ~ \frac{1}{2}\sum_{i=1}^m y_i+\frac{1}{2} \\
        = & ~ 1,
    \end{align*}
    therefore, we can apply Lemma~\ref{lem:scp_mw_bound} for the following bound, let $p\in {\cal K}_s$, then
    \begin{align*}
        \frac{1}{2}-\frac{\alpha}{6\rho} \leq & ~ \frac{1}{T}\sum_{t=1}^T \langle \ell^t, p\rangle+\eta+\frac{\log m}{\eta T} \\
        = & ~ \frac{1}{T}\sum_{t=1}^T p_i\cdot (\frac{1}{2\rho}(b_i-\langle a_i,x^t\rangle)+\frac{1}{2})+\eta + \frac{\log m}{\eta T},
    \end{align*}
    rearranging gives
    \begin{align*}
        -\frac{\alpha}{6\rho} \leq & ~ \frac{1}{T}\sum_{t=1}^T \sum_{i=1}^m p_i\cdot \frac{1}{2\rho} (b_i-\langle a_i,x^t\rangle)+\eta+\frac{\log m}{\eta T},
    \end{align*}
    recall that we set the final output as $\ov x=1/T\sum_{t=1}^T x^t$, the above bound could further be written as
    \begin{align*}
        \sum_{i=1}^m p_i\cdot \langle a_i, \ov x\rangle\leq & ~ \sum_{i=1}^m p_i\cdot b_i+2\rho\eta+\frac{2\rho\log m}{\eta T}+\frac{\alpha}{3},
    \end{align*}
    picking $\eta=\sqrt{\frac{\log m}{T}}=\frac{\alpha}{6\rho}$ and $T=\frac{36\rho^2\log m}{\alpha^2}$ (so that $\eta\leq 1/2$, as $\alpha\leq 3\rho$) gives
    \begin{align*}
        \sum_{i=1}^m p_i\cdot \langle a_i, \ov x\rangle\leq & ~ \sum_{i=1}^m p_i\cdot b_i+\alpha.
    \end{align*}
    As the above bound holds for any $p\in {\cal K}_s$, it must be the case that $\ov x$ satisfies all but at most $s-1$ constraints with additive error $\alpha$. Suppose otherwise, there are $s$ constraints violate the $\alpha$-additive error condition, then we could set $p$ as the uniform distribution over these $s$ constraints and this implies that these $s$ constraints satisfy the $\alpha$-additive error condition, a contradiction.
\end{proof}

\subsection{Privacy Guarantee of Dense Multiplicative Weights Update}\label{sec:achieve_constraint_private}
Next, we prove that Algorithm~\ref{alg:constraint_mwu} is private as long as the oracle is private. The algorithm only utilizes the private data via the oracle minimization, hence, as long as the oracle is privately minimizes over $\K$ at each step $t\in [T]$, the final output will automatically be private as $\K$ is convex. To do so, we recall a lemma first proved in~\cite{hru13}, showing that for two neighboring databases that differ by one action, then their projected distributions $\wt y$ and $\wt y'$ satisfy $\|\wt y-\wt y'\|_1\leq 2/s$.

\begin{lemma}[\cite{hru13}]\label{lem:pro_dis_bound}
     Let $F: \mathcal{F} \rightarrow [0,1]$ and $F':\mathcal{F} \cup \{f'\} \rightarrow [0,1]$ be two measures over their respective set of actions. Let density parameter $s \in \mathbb{N}$ and it satisfies that (1) $|F|,|F'| \leq s$, (2) $F_f = F'_f$ for every $f \in \mathcal{F}$. Let $\wt{F},\wt{F}'$ be the corresponding Bregman projections onto the set of $1/s$-dense distributions, then we have
     \begin{align*}
         \|\wt{F} - \wt{F}'\|_1 \leq 2/s
     \end{align*}
\end{lemma}

Now we are able to present our main theorem for high sensitivity constraint privacy:
\begin{theorem}[Privacy guarantee]\label{thm:constraint_privacy_guarantee}
Let $\epsilon\in(0,1]$ and $\delta\in(0,1/2]$. Consider a symmetric cone program with constraint matrix
$A\in\mathcal{J}^m$, vector $b\in\mathbb{R}^m$, and width $\rho$. Fix parameters
$\alpha\in(0,3\rho]$, $\beta\in(0,1)$, and $T= 36\alpha^{-2}\rho^2\log m$. Given access to an $(\alpha/3,\beta/T)$-approximate, $\rho$-bounded oracle that is also
$\epsilon'$-private where $\epsilon'=\frac{\epsilon}{\sqrt{8T\log(1/\delta)}}$. For any neighboring instances, the oracle inputs satisfy
\begin{align*}
    \|\widetilde{y}\|_\infty \leq \tfrac{1}{s}, \qquad
    \|\widetilde{y}'\|_\infty \leq \tfrac{1}{s}, \qquad
    \|\widetilde{y}-\widetilde{y}'\|_1 \le \tfrac{2}{s}.
\end{align*}
Then Algorithm~\ref{alg:constraint_mwu}, run with density parameter $s$, is
$(\epsilon,\delta)$-high sensitivity constraint private.
\end{theorem}
\begin{proof}
Since the oracle is $\epsilon'$-private, the overall $(\epsilon,\delta)$-constraint private is achieved via adaptive composition (Lemma~\ref{lem:composition_lemma}). When adding or removing a constraint, we note that $A$ and $A'$ are identical except for one additional constraint. Since we project the distribution onto $\K_s$, it must be the case $\|\wt y\|_\infty,\|\wt y'\|_\infty\leq 1/s$, and $\|\wt y-\wt y'\|_1\leq 2/s$ follows from Lemma~\ref{lem:pro_dis_bound}. Thus, since two distributions are identical except for the probability associated with the additional constraint, we have completed the proof.
\end{proof}

\begin{remark}[Data-independent number of rounds]\label{rem:public_rounds}
The number of rounds $T$ and the step size $\eta$ of Algorithm~\ref{alg:constraint_mwu} depend on $\log m$. Neighboring instances, however, have $m$ and $m+1$ constraints. The proof of Theorem~\ref{thm:constraint_privacy_guarantee} needs $T$ and $\eta$ to be the same on both instances, so that the two runs apply the same updates to the common constraints. We therefore compute $T$ and $\eta$ from a public upper bound $\bar m\geq m+1$ on the number of constraints. Lemma~\ref{lem:approximate_feasibility} still holds with $\log m$ replaced by $\log\bar m$, because its proof uses $\log m$ only to bound the term $\log|{\cal F}|/(\eta T)$ of Lemma~\ref{lem:scp_mw_bound}, and $\log|{\cal F}|\leq\log\bar m$. In Theorems~\ref{thm:private_covering_sdp} and~\ref{thm:private_covering_scp}, $m$ in $T$, in $\eta$ and in the condition on $s$ stands for $\bar m$.
\end{remark}

\subsection{Application I: Private Solver for Covering SDP}

As an application, we show how to develop a private solver for covering semidefinite programming:
\begin{align*}
    \min_{X\in \R^{r\times r}} & ~ \Tr(X) \\
    \text{s.t.} & ~ \langle A_i, X\rangle\geq 1, \forall i\in [m] \\
    & ~ X\succeq 0
\end{align*}
where $A_1,\ldots,A_m\succeq 0$. Again, we consider the convex set $\K_{\OPT}=\{X\succeq 0: \Tr(X)=\OPT \}$, and we assume that $\OPT$ is given to the algorithm as a public input, the same on neighboring instances. Remark~\ref{rem:public_opt} explains how a public upper bound on $\OPT$ can be used instead. We also assume that every admissible constraint matrix, including the extra constraint of a neighboring instance, satisfies $A\succeq 0$ and $\|A\|\leq 1$. A public bound $B$ on the norms $\|A_i\|$ reduces to the case $B=1$ by replacing $A_i$ with $A_i/B$, which multiplies $\OPT$ by $B$. A scale computed from the private $A_i$ would not be private. Since $A_1\preceq I$, every feasible $X$ satisfies $\Tr(X)\geq\langle A_1,X\rangle\geq 1$, so $\OPT\geq 1$. Thus, we can consider the feasibility problem
\begin{align*}
    \text{find} & ~ X\in \K_{\OPT} \\
    \text{s.t.} & ~ \langle A_i, X\rangle\geq 1, \forall i\in [m].
\end{align*}
It is natural to study high sensitivity constraint privacy in this setting, as one could view $A_i$ as the positive semidefinite constraint matrix associated with each individual data point as in the case of robust mean estimation~\cite{cdg19}, robust covariance estimation~\cite{cdgw19} and E-optimal experimental design~\cite{vg19}. 

We need to design an oracle to solve the maximization problem
\begin{align*}
    O(y) = & ~ \arg\max_{X\in \K_{\OPT}} \langle \sum_{i=1}^m y_i A_i, X\rangle,
\end{align*}
this is a linear maximization problem, and $\K_{\OPT}$ is the intersection of a hyperplane and the positive semidefinite cone, and the solutions to the linear maximization are the extreme rays of the cone:
\begin{align*}
    X^* = & ~ uu^\top,
\end{align*}
where $\|u\|_2^2=\OPT$. To implement a private oracle, it is tempting to use exponential mechanism directly, but note that there are infinitely many $r$-dimensional vectors $u$ satisfying $\|u\|_2^2=\OPT$, so the accuracy guarantee of Lemma~\ref{lem:acc_gua_exp_mec} becomes meaningless as $|{\cal R}|=\infty$. To address this issue, we use a $\gamma$-net argument to quantize the sphere $S=\{u\in \R^r: \|u\|_2^2= \OPT \}$ into finitely many points. 

\begin{lemma}[\cite{lt91,v18}]
\label{lem:net_argument}
Let $B=\{u\in \R^r: \|u\|_2\leq R\}$, there exists a finite collection of points $N\subset \R^r$ such that for any $u\in B$, there exists $v\in N$ such that $\|u-v\|_2\leq \gamma$ for $\gamma>0$. Moreover, one can take $N\subseteq B$ with $|N|\leq(1+2R/\gamma)^r$. The same bound holds for a net $N\subseteq S$ of the sphere $S=\{u\in \R^r: \|u\|_2=R\}$, by taking $N$ to be a maximal $\gamma$-separated subset of $S$.
\end{lemma}
Our private oracle will then be the exponential mechanism, over the net $N$.

\begin{lemma}\label{lem:private_oracle}
Let $\OPT,\epsilon>0$, $\beta\in(0,1)$, $s\in\mathbb{N}$, and assume that every admissible constraint matrix, including the extra constraint of a neighboring instance, satisfies $A\succeq 0$ and $\|A\|\leq 1$. Assume neighboring inputs $y,y' \in \mathbb{R}^m$ satisfy $\|y\|_\infty \leq \frac{1}{s}$, $\|y'\|_\infty \leq \frac{1}{s}$ and $\|y-y'\|_1 \leq \frac{2}{s}$. Let $\alpha_0\in(0,\OPT]$, let $S=\{u\in\mathbb{R}^r:\|u\|_2^2=\OPT\}$, and let $N\subseteq S$ be a $\gamma$-net of $S$ with $\gamma=\alpha_0/(2\sqrt{\OPT})$ and $|N|\leq(1+2\sqrt{\OPT}/\gamma)^r\leq(5\OPT/\alpha_0)^r$ (a maximal $\gamma$-separated subset of $S$). Let
$O(y)$ denote the exponential mechanism on $N$ with privacy parameter $\epsilon$ and sensitivity parameter $\Delta=3\OPT/s$ that outputs $uu^\top$ for the selected $u$, with quality score
\begin{align*}
  Q(u,y) = & ~ \langle \sum_{i=1}^m y_i A_i, uu^\top \rangle - 1 .
\end{align*}
Then $O$ is $\epsilon$-private with respect to the neighboring inputs above (constraint sets that differ by one extra admissible constraint, and $y,y'$ as assumed), every output of $O(y)$ lies in $\K_{\OPT}$, and $O(y)$ is an $(\alpha,\beta)$-approximate, $\rho$-bounded oracle with
\begin{align*}
  \alpha = \alpha_0+\frac{6\OPT}{s \epsilon}\Big(r\log\frac{5\OPT}{\alpha_0}+\log\frac{1}{\beta}\Big), \qquad
  \rho \le \max\{1,\OPT\} .
\end{align*}
\end{lemma}

\begin{proof}
Every $u\in N$ has $\|u\|_2^2=\OPT$, so $uu^\top\succeq 0$ and $\Tr(uu^\top)=\OPT$, that is, every output lies in $\K_{\OPT}$. We next bound the width. Since $A_i\succeq 0$ and $\|A_i\|\leq 1$, every $u\in N$ satisfies
\begin{align*}
    0\leq \langle A_i,uu^\top\rangle = & ~ u^\top A_i u \leq \|u\|_2^2 = \OPT,
\end{align*}
hence $|\langle A_i,uu^\top\rangle-1|\leq \max\{1,\OPT\}$. The same bound holds for every $X\in\K_{\OPT}$, since $0\leq\langle A_i,X\rangle\leq\|A_i\|\Tr(X)\leq\OPT$.

For accuracy guarantee, we first need to compute the sensitivity of the quality score. To do so, note that for the first $m$ entries, we have $\sum_{i=1}^m |y_i-y_i'|\leq 2/s$ and $D'$ might have one more constraint, therefore the extra entry of $y'$ has its magnitude being at most $1/s$. We let $y, y'\in \R^{m+1}$ and set $y_{m+1}=0$, the sensitivity is
\begin{align*}
    |\langle \sum_iy_iA_i-\sum_i y_i'A_i, uu^\top\rangle | = & ~ |u^\top (\sum_iy_iA_i-\sum_i y_i'A_i)u | \\
    = & ~ | \sum_i (y_i-y_i') u^\top A_i u|\\
    \leq & ~ \max_{i} |u^\top A_iu| \cdot \|y-y'\|_1 \\
    \leq & ~ \max_i \|A_i\|\cdot \|u\|_2^2\cdot 3/s \\
    \leq & ~ \frac{3\OPT}{s}.
\end{align*}
So the quality score has sensitivity at most $\Delta=3\OPT/s$. As $\OPT$ and $s$ are public, $\Delta$ and the net $N$ do not depend on the private data. Hence the exponential mechanism is $\epsilon$-private by Definition~\ref{def:exponential_mechanism}. Next we bound the error of the discretization. Let $M=\sum_{i=1}^m y_iA_i$. As $y$ is a distribution, $\|M\|\leq \sum_{i=1}^m y_i\|A_i\|\leq 1$. The maximum of $\langle M,X\rangle$ over $\K_{\OPT}$ is attained at an extreme point $X=vv^\top$ with $v\in S$. Pick $u\in N$ with $\|u-v\|_2\leq\gamma$. As $M$ is symmetric,
\begin{align*}
    |\langle M,vv^\top\rangle-\langle M,uu^\top\rangle| = & ~ |(v-u)^\top M(v+u)| \\
    \leq & ~ \|v-u\|_2\cdot\|M\|\cdot\|v+u\|_2 \\
    \leq & ~ 2\gamma\sqrt{\OPT} = \alpha_0 ,
\end{align*}
so $\max_{u\in N}Q(u,y)\geq \max_{X\in\K_{\OPT}}\langle M,X\rangle-1-\alpha_0$. Finally, by Lemma~\ref{lem:acc_gua_exp_mec} with sensitivity $3\OPT/s$ and $|N|\leq(5\OPT/\alpha_0)^r$, with probability at least $1-\beta$ the selected $u$ satisfies
\begin{align*}
    Q(u,y) \geq & ~ \max_{u'\in N}Q(u',y)-\frac{6\OPT}{s\cdot \epsilon}\log(|N|/\beta) \\
    \geq & ~ \max_{X\in\K_{\OPT}}\langle M,X\rangle-1-\alpha_0-\frac{6\OPT}{s\cdot \epsilon}\Big(r\log\frac{5\OPT}{\alpha_0}+\log\frac{1}{\beta}\Big) \\
    = & ~ \max_{X\in\K_{\OPT}}\langle M,X\rangle-1-\alpha .
\end{align*}
This is the accuracy required of the oracle, since the covering constraints read $\langle -A_i,X\rangle\leq -1$ in the form used by Algorithm~\ref{alg:constraint_mwu}. We hence complete the proof.
\end{proof}

We are set in a position to state the final privacy and utility guarantee of our algorithm for high sensitivity constraint privacy. 

\begin{theorem}[Formal version of Theorem~\ref{thm:private_convering_sdp_informal}]
\label{thm:private_covering_sdp}
Consider a feasible covering SDP with $m\geq 3$ constraints, where every admissible constraint matrix satisfies $A\succeq 0$ and $\|A\|\leq 1$, and assume that its optimal value $\OPT$ is given to the algorithm as a public input, the same on neighboring instances. Let $\epsilon\in(0,1]$, $\delta\in(0,1/2]$, $\beta\in (0,1)$ and $\alpha\in(0,\OPT]$. Run Algorithm~\ref{alg:constraint_mwu} on $\K_{\OPT}=\{X\succeq 0:\Tr(X)=\OPT\}$ with width $\rho=\OPT$, accuracy $\alpha$, $T=\lceil 36\alpha^{-2}\OPT^2\log m\rceil$ rounds (with $m$ read as the public bound $\bar m$ of Remark~\ref{rem:public_rounds}, here and in the condition on $s$) and density parameter $s$, using the oracle $O(y)$ of Lemma~\ref{lem:private_oracle} with $\alpha_0=\alpha/6$, privacy parameter $\epsilon'=\epsilon/\sqrt{8T\log(1/\delta)}$ and failure probability $\beta/T$. We also assume that every instance has at least $s$ constraints, so that the $1/s$-dense distributions used by Algorithm~\ref{alg:constraint_mwu} exist. There is an absolute constant $C_0$ such that, if
\begin{align*}
    s \geq & ~ C_0\cdot\frac{\OPT^2}{\epsilon\,\alpha^2}\cdot \log^{1/2}(1/\delta)\log^{1/2} m\cdot\Big(r\log\frac{2\OPT}{\alpha}+\log\frac{\log m}{\beta}\Big),
\end{align*}
then, for every fixed value of this input that is at least $1$ (Remark~\ref{rem:public_opt}), Algorithm~\ref{alg:constraint_mwu} is $(\epsilon,\delta)$-private with respect to high sensitivity constraint privacy, its output $X^*$ satisfies $X^*\succeq 0$ and $\Tr(X^*)=\OPT$, and with probability at least $1-\beta$ fewer than $s$ constraints have $\langle A_i, X^*\rangle<1-\alpha$.
\end{theorem}

\begin{proof}
Recall that $\OPT\geq 1$. Every $X\in\K_{\OPT}$ and every admissible constraint matrix $A$ satisfy $0\leq\langle A,X\rangle\leq\|A\|\Tr(X)\leq\OPT$, so $|\langle A,X\rangle-1|\leq\max\{1,\OPT\}=\OPT$. Hence $\rho=\OPT$ is a valid width, and it is public. An optimal solution lies in $\K_{\OPT}$, so the feasibility program is feasible. As $\alpha\leq\OPT=\rho$, we have $\alpha\leq 3\rho$. To apply Lemma~\ref{lem:approximate_feasibility} and Theorem~\ref{thm:constraint_privacy_guarantee}, we need an $(\alpha/3,\beta/T)$-approximate oracle that is $\epsilon'$-private, where $T=36\alpha^{-2}\rho^2\log m$ and $\epsilon'=\frac{\epsilon}{\sqrt{8T\log(1/\delta)}}$. We use the oracle of Lemma~\ref{lem:private_oracle} with $\alpha_0=\alpha/6$, privacy parameter $\epsilon'$ and failure probability $\beta/T$. Its error is at most $\alpha/3$ provided
\begin{align*}
    \frac{6\OPT}{s\cdot \epsilon'}\Big(r\log\frac{30\OPT}{\alpha}+\log\frac{T}{\beta}\Big) \leq & ~ \frac{\alpha}{6}.
\end{align*}
As $\frac{1}{\epsilon'}=\frac{\sqrt{8T\log(1/\delta)}}{\epsilon}$ and $\sqrt{T}=\frac{6\OPT\sqrt{\log m}}{\alpha}$, this condition is equivalent to
\begin{align*}
    s \geq & ~ \frac{216\sqrt{8}\,\OPT^2}{\alpha^2\epsilon}\cdot \sqrt{\log(1/\delta)\log m}\cdot \Big(r\log\frac{30\OPT}{\alpha}+\log\frac{36\OPT^2\log m}{\alpha^2\beta}\Big).
\end{align*}
We bound the last factor. As $\alpha\leq\OPT$, we have $\log(2\OPT/\alpha)\geq\log 2$, hence $\log 15\leq 4\log(2\OPT/\alpha)$ and $\log 36\leq 6\log(2\OPT/\alpha)$. As $m\geq 3$ and $\beta<1$, we have $\log(\log m/\beta)>0$. Using also $r\geq 1$ and $\log(\OPT/\alpha)\leq\log(2\OPT/\alpha)$, we get
\begin{align*}
    r\log\frac{30\OPT}{\alpha} = & ~ r\Big(\log 15+\log\frac{2\OPT}{\alpha}\Big) \leq 5r\log\frac{2\OPT}{\alpha}, \\
    \log\frac{36\OPT^2\log m}{\alpha^2\beta} = & ~ \log 36+2\log\frac{\OPT}{\alpha}+\log\frac{\log m}{\beta} \leq 8r\log\frac{2\OPT}{\alpha}+\log\frac{\log m}{\beta}.
\end{align*}
So the last factor is at most $13\,(r\log(2\OPT/\alpha)+\log(\log m/\beta))$, and the condition on $s$ in the theorem with $C_0=13\cdot 216\sqrt{8}$ implies the condition above. So far we treated $T$ as a real number, but the algorithm runs $\lceil T\rceil$ rounds and computes $\epsilon'$ from $\lceil T\rceil$. Rounding $T$ up keeps Lemma~\ref{lem:approximate_feasibility} valid, since then $\eta=\sqrt{\log m/T}\leq\alpha/(6\rho)$. Rounding up at most doubles $T$, as $36\OPT^2\log m/\alpha^2\geq 36\log 3>1$, so it changes only the absolute constant $C_0$. Under the condition of the theorem, the oracle is therefore $\epsilon'$-private, $(\alpha/3,\beta/T)$-approximate and $\OPT$-bounded, and its outputs lie in $\K_{\OPT}$. Lemma~\ref{lem:approximate_feasibility}, applied to the constraints $\langle -A_i,X\rangle\leq -1$, then shows that with probability at least $1-\beta$ fewer than $s$ constraints have $\langle A_i,X^*\rangle<1-\alpha$. The output $X^*$ is an average of oracle outputs, so it lies in the convex set $\K_{\OPT}$, that is, $X^*\succeq 0$ and $\Tr(X^*)=\OPT$. Finally, Theorem~\ref{thm:constraint_privacy_guarantee}, which uses $\epsilon\leq 1$ and $\delta\leq 1/2$, gives $(\epsilon,\delta)$-privacy. This completes the proof.
\end{proof}

\begin{remark}[Public value of $\OPT$]\label{rem:public_opt}
Theorems~\ref{thm:private_covering_sdp} and~\ref{thm:private_covering_scp} take $\OPT$ as public input. Privacy holds for every value of this input that is at least $1$, because the privacy proofs use only that this value is the same on neighboring instances and at least $1$. This lower bound keeps the losses of Algorithm~\ref{alg:constraint_mwu} in $[0,1]$. Only the accuracy guarantee needs the value to equal the optimal value of the given instance, or, as explained below, to be at least that value. A binary search for $\OPT$ would read the private constraints, and our privacy analysis does not cover it. If only a public upper bound $V\geq\OPT$ is known, both theorems hold with $\OPT$ replaced by $V$ in the algorithm, in the trace of the output and in the condition on $s$, for every $\alpha\in(0,V]$. To see this, note first that if $X_{\rm opt}$ is an optimal solution, then $(V/\OPT)X_{\rm opt}\in\K_V$ satisfies every constraint, since $\langle A_i,X_{\rm opt}\rangle\geq 1$ and $V/\OPT\geq 1$. The same holds for covering SCPs. The proofs use $\OPT$ only through the feasibility of $\K_{\OPT}$, through $\OPT\geq 1$ and through Lemmas~\ref{lem:private_oracle} and~\ref{lem:private_oracle_eja}, which hold for any positive value in place of $\OPT$.
\end{remark}

\subsection{Application II: Private Solver for Covering SCP}

More generally, consider the covering symmetric cone programming:
\begin{align*}
    \min_{x\in \J} & ~ \Tr(x) \\
    \text{s.t.} & ~ \langle a_i, x\rangle\geq 1, \forall i\in [m] \\
    & ~ x\in \K,
\end{align*}
where $a_1,\ldots,a_m\in \K$. As in Application I, we assume that $\OPT$ is given to the algorithm as a public input, the same on neighboring instances. We also assume that every admissible constraint, including the extra constraint of a neighboring instance, satisfies $a\in\K$ and $\|a\|_\infty\leq 1$. Every feasible $x$ satisfies $1\leq\langle a_1,x\rangle\leq\|a_1\|_\infty\|x\|_1\leq\|x\|_1=\Tr(x)$ by Lemma~\ref{lem:holder}, so $\OPT\geq 1$. Setting $\K_{\OPT}=\{x\in \K: \Tr(x)=\OPT \}$, we reduce the optimization to feasibility:
\begin{align*}
    \text{find} & ~ x\in \K_{\OPT} \\
    \text{s.t.} & ~ \langle a_i, x\rangle\geq 1, \forall i\in [m].
\end{align*}
Recall that we need to design an oracle to solve
\begin{align*}
    O(y) = & ~ \arg\max_{x\in \K_{\OPT}} \langle\sum_{i=1}^m y_ia_i, x\rangle,
\end{align*}
in the case of SDP, we noted that the optimal solutions to a linear maximization over the positive semidefinite cone are the extreme rays of the cone, i.e., rank-1 positive semidefinite matrices. For SCP, we consider the case where $\J$ is a simple Jordan algebra, which covers the interesting cases including real symmetric matrices, Hermitian symmetric matrices and spin factors (the algebra whose cone of the squares is the second-order cone). In this scenario, it's easy to see that the optimal solutions are $\OPT\cdot q$ for a primitive idempotent $q$. To use a $\gamma$-net argument, we need to understand the dimension of all the primitive idempotents in $\J$. These primitive idempotents form a connected and compact manifold, whose dimension can be characterized as follows:
\begin{lemma}[\cite{fk94}]
Let $\J$ be a simple Jordan algebra with rank $r$ and Peirce constant $d$. Let $Q=\{ q\in \J: \text{$q$ is a primitive idempotent}\}$, then $\dim(Q)=d(r-1)$. Let $C=\{c\cdot q: c\in \R, q\in Q \}$, then $\dim(C)=d(r-1)+1$.
\end{lemma}
Peirce constant is a parameter related to the Pierce decomposition of a Jordan frame. We list $r$ and $d$ for all the simple Jordan algebras.
\begin{table}[!ht]
    \centering
    \begin{tabular}{|l|l|l|l|}
    \hline
      Algebra $\J$ & Rank $r$ & Peirce constant $d$ & $k_C=d(r-1)+1$ \\ \hline
       $r\times r$ real symmetric matrices & $r$ & 1 & $r$ \\ \hline
       $r\times r$ complex Hermitian matrices & $r$ & 2 & $2r-1$ \\ \hline
       $r\times r$ quaternionic Hermitian matrices & $r$ & 4 & $4r-3$ \\ \hline
       spin factor $\R\times\R^n$, $n\geq 2$ & 2 & $n-1$ & $n$ \\ \hline
       Albert algebra & 3 & 8 & 17 \\ \hline
    \end{tabular}
    \vspace{2mm}
    \caption{Five types of simple Jordan algebras, their rank $r$, their Peirce constant $d$, and the dimension $k_C$ of the cone $C$ of scaled primitive idempotents.}
    \label{tab:simple_jordan}
\end{table}

The spin factor $\R\times\R^n$ consists of the vectors $(x_0,\bar x)\in\R^{n+1}$, and its cone of squares is the second-order cone. We can hence conclude that for all simple Jordan algebras, the dimension of $C$ is $k_C=d(r-1)+1$, which is at most $4r$ for the three matrix algebras, equals $n$ for the spin factor $\R\times\R^n$ (whose rank is $2$), and equals $17$ for the Albert algebra (rank $3$). Lemma~\ref{lem:private_oracle_eja} below is the analog of Lemma~\ref{lem:private_oracle}. Its net has size $\exp(O(k_C\log(2\OPT/\alpha_0)))$, with no further dependence on the rank. For the matrix algebras, this bound uses that the difference of two primitive idempotents has at most two nonzero eigenvalues. For the spin factors and the Albert algebra, it uses that their rank is at most $3$. The proof of Theorem~\ref{thm:private_covering_sdp}, run with Lemma~\ref{lem:private_oracle_eja} in place of Lemma~\ref{lem:private_oracle}, gives Theorem~\ref{thm:private_covering_scp}.

\begin{lemma}\label{lem:private_oracle_eja}
Let $\J$ be a simple EJA of rank $r$ and Peirce constant $d$, and assume that every admissible constraint, including the extra constraint of a neighboring instance, satisfies $a\in\K$ and $\|a\|_\infty\leq 1$. Let $\OPT,\epsilon>0$, $\beta\in (0,1)$, $s\in \mathbb{N}$, $\alpha_0\in(0,\OPT]$, let $k_C=d(r-1)+1$, and let $Q_{\OPT}=\{\OPT\cdot q: q\in Q \}\subseteq \K_{\OPT}$. Assume that neighboring inputs $y,y'\in \R^m$ satisfy $\|y\|_\infty, \|y'\|_\infty\leq 1/s$ and $\|y-y'\|_1\leq 2/s$. Then there is a set $N\subseteq Q_{\OPT}$, depending only on $\J$, $\OPT$ and $\alpha_0$, with $\log|N|\leq 2k_C\log(7\OPT/\alpha_0)$, such that every $q\in Q$ has some $p\in N$ with $\|p-\OPT\cdot q\|_1\leq\alpha_0$. Let $O(y)$ denote the exponential mechanism on $N$ with privacy parameter $\epsilon$, sensitivity parameter $\Delta=3\OPT/s$ and quality score
\begin{align*}
    Q(p, y)= & ~ \langle \sum_{i=1}^m y_ia_i, p\rangle-1.
\end{align*}
Then $O$ is $\epsilon$-private with respect to the neighboring inputs above (constraint sets that differ by one extra admissible constraint, and $y,y'$ as assumed), every output of $O(y)$ lies in $\K_{\OPT}$, and $O(y)$ is an $(\alpha,\beta)$-approximate, $\rho$-bounded oracle with
\begin{align*}
    \alpha =  \alpha_0+\frac{6\OPT}{s\cdot \epsilon}\Big(2k_C\log\frac{7\OPT}{\alpha_0}+\log\frac{1}{\beta}\Big), & \quad
    \rho \leq  \max\{1,\OPT\}.
\end{align*}
\end{lemma}

\begin{proof}
We first construct $N$. Here $\|x\|_1$ and $\|x\|_2$ are the norms of Definition~\ref{def:inner_product}. A Jordan algebra isomorphism preserves eigenvalues, hence these norms, the trace and the set $Q$. If $r=1$, then $\J=\R$, $Q=\{e\}$, and $N=Q_{\OPT}$ is a single point. Let $r\geq 2$. By the classification of simple Euclidean Jordan algebras~\cite{fk94}, $\J$ is isomorphic to one of the algebras of Table~\ref{tab:simple_jordan}. We treat the three kinds in turn.

\emph{Matrix algebras.} Let $\J$ be the algebra of $r\times r$ Hermitian matrices over $\mathbb{F}\in\{\R,\mathbb{C},\mathbb{H}\}$, so that $d=\dim_{\R}\mathbb{F}\in\{1,2,4\}$ and $\mathbb{F}^r\cong\R^{dr}$. Then $Q=\{uu^*: u\in\mathbb{F}^r,\ \|u\|_2=1\}$, the eigenvalues of an element are its eigenvalues as a Hermitian matrix, and $\|x\|_2$ is the Frobenius norm. For unit vectors $u,v\in\mathbb{F}^r$, the Hermitian matrix $X=uu^*-vv^*$ maps $\mathbb{F}^r$ into the right $\mathbb{F}$-span of $u$ and $v$. Write $X=\sum_j\lambda_jw_jw_j^*$ with orthonormal $w_j$ and real $\lambda_j$. Every $w_j$ with $\lambda_j\neq 0$ equals $X(w_j\lambda_j^{-1})$, so it lies in the span of $u$ and $v$. Since orthonormal vectors are linearly independent, at most two $\lambda_j$ are nonzero. As $\Tr(X)=\|u\|_2^2-\|v\|_2^2=0$, they are $\lambda$ and $-\lambda$ for some $\lambda\geq 0$. As $uu^*-vv^*=u(u-v)^*+(u-v)v^*$,
\begin{align*}
    \|uu^*-vv^*\|_1 = \sqrt{2}\,\|uu^*-vv^*\|_2 \leq \sqrt{2}\,\big(\|u(u-v)^*\|_2+\|(u-v)v^*\|_2\big) = 2\sqrt{2}\,\|u-v\|_2 .
\end{align*}
Let $N_0$ be a maximal $\gamma$-separated subset of the unit sphere of $\R^{dr}$ with $\gamma=\alpha_0/(2\sqrt{2}\,\OPT)$, and let $N=\{\OPT\cdot vv^*: v\in N_0\}$. By Lemma~\ref{lem:net_argument}, every unit vector $u$ has some $v\in N_0$ with $\|u-v\|_2\leq\gamma$, so that $\|\OPT\cdot uu^*-\OPT\cdot vv^*\|_1\leq 2\sqrt{2}\,\OPT\gamma=\alpha_0$. By the same lemma, the net also satisfies $|N|\leq(1+4\sqrt{2}\,\OPT/\alpha_0)^{dr}\leq(7\OPT/\alpha_0)^{2k_C}$, because $\alpha_0\leq\OPT$ and $dr\leq 2(d(r-1)+1)$ for $r\geq 2$.

\emph{Spin factors.} Let $\J=\R\times\R^n$ with $n\geq 2$, so that $r=2$, $d=n-1$ and $k_C=n$. The element $(x_0,\bar x)$ has eigenvalues $x_0\pm\|\bar x\|_2$, and $Q=\{\frac{1}{2}(1,w): w\in\R^n,\ \|w\|_2=1\}$. For $q=\frac{1}{2}(1,w)$ and $q'=\frac{1}{2}(1,w')$, the difference $q-q'=(0,\frac{1}{2}(w-w'))$ has eigenvalues $\pm\frac{1}{2}\|w-w'\|_2$, so $\|q-q'\|_1=\|w-w'\|_2$. Let $N_0$ be a maximal $(\alpha_0/\OPT)$-separated subset of the unit sphere of $\R^n$, and let $N=\{\frac{\OPT}{2}(1,w): w\in N_0\}$. By Lemma~\ref{lem:net_argument}, every $q\in Q$ has some $p\in N$ with $\|p-\OPT\cdot q\|_1\leq\alpha_0$. By the same lemma, the net also satisfies $|N|\leq(1+2\OPT/\alpha_0)^{n}\leq(7\OPT/\alpha_0)^{2k_C}$.

\emph{The Albert algebra.} Here $r=3$, $d=8$, $k_C=17$ and $\dim(\J)=27$. Every element has three eigenvalues, so $\|x\|_1\leq\sqrt{3}\,\|x\|_2$ by the Cauchy--Schwarz inequality. Also, $\|\OPT\cdot q\|_2=\OPT$ for $q\in Q$. Let $N$ be a maximal $(\alpha_0/\sqrt{3})$-separated subset of $Q_{\OPT}$ in the norm $\|\cdot\|_2$. By maximality, every $\OPT\cdot q$ has some $p\in N$ with $\|p-\OPT\cdot q\|_2\leq\alpha_0/\sqrt{3}$, and hence $\|p-\OPT\cdot q\|_1\leq\alpha_0$. Through an isometry $\phi:\J\to\R^{27}$ (Definition~\ref{def:dim_iso}), the balls of radius $\alpha_0/(2\sqrt{3})$ around the points of $N$ are disjoint and lie in the ball of radius $\OPT+\alpha_0/(2\sqrt{3})$, so $|N|\leq(1+2\sqrt{3}\,\OPT/\alpha_0)^{27}\leq(7\OPT/\alpha_0)^{2k_C}$.

In all cases $\log|N|\leq 2k_C\log(7\OPT/\alpha_0)$, and $N$ depends only on $\J$, $\OPT$ and $\alpha_0$.

Every output lies in $Q_{\OPT}\subseteq\K_{\OPT}$. For $x\in\K_{\OPT}$ and an admissible constraint $a\in\K$ with $\|a\|_\infty\leq 1$, self-duality of $\K$ and Lemma~\ref{lem:holder} give $0\leq\langle a,x\rangle\leq\|a\|_\infty\|x\|_1\leq\|x\|_1=\Tr(x)=\OPT$, where $\|x\|_1=\Tr(x)$ because $x\in\K$ has nonnegative eigenvalues. Hence $|\langle a,x\rangle-1|\leq\max\{1,\OPT\}$. As in the proof of Lemma~\ref{lem:private_oracle}, it follows that the quality score has sensitivity at most $\Delta=3\OPT/s$. As $\OPT$ and $s$ are public and $N$ depends only on $\J$, $\OPT$ and $\alpha_0$, the exponential mechanism is $\epsilon$-private by Definition~\ref{def:exponential_mechanism}.

Let $M=\sum_{i=1}^m y_ia_i$, so that $\|M\|_\infty\leq 1$. Every $x\in\K_{\OPT}$ has a spectral decomposition $x=\sum_j\lambda_jq_j$ with $\lambda_j\geq 0$ and $\sum_j\lambda_j=\OPT$, so it is a convex combination of the points $\OPT\cdot q_j\in Q_{\OPT}$. Hence the maximum of $\langle M,x\rangle$ over the compact set $\K_{\OPT}$ is attained at a point $\OPT\cdot q^*$ with $q^*\in Q$. Pick $p\in N$ with $\|p-\OPT\cdot q^*\|_1\leq\alpha_0$. By Lemma~\ref{lem:holder},
\begin{align*}
    |\langle M,p-\OPT\cdot q^*\rangle| \leq & ~ \|M\|_\infty\cdot\|p-\OPT\cdot q^*\|_1 \leq \alpha_0 .
\end{align*}
So $\max_{p'\in N}Q(p',y)\geq\max_{x\in\K_{\OPT}}\langle M,x\rangle-1-\alpha_0$. Lemma~\ref{lem:acc_gua_exp_mec} with sensitivity $\Delta=3\OPT/s$ then shows that, with probability at least $1-\beta$, the selected $\hat p$ satisfies
\begin{align*}
    Q(\hat p,y) \geq & ~ \max_{x\in\K_{\OPT}}\langle M,x\rangle-1-\alpha_0-\frac{6\OPT}{s\cdot\epsilon}\log(|N|/\beta) \geq \max_{x\in\K_{\OPT}}\langle M,x\rangle-1-\alpha,
\end{align*}
using $\log|N|\leq 2k_C\log(7\OPT/\alpha_0)$. As in the proof of Lemma~\ref{lem:private_oracle}, this is the accuracy required of the oracle.
\end{proof}

\begin{theorem}[Formal version of Theorem~\ref{thm:private_convering_scp_informal}]
\label{thm:private_covering_scp}
Let $\J$ be a simple EJA of rank $r$ and Peirce constant $d$, and let $k_C=d(r-1)+1$. Consider a feasible covering SCP over $\J$ with $m\geq 3$ constraints, where every admissible constraint satisfies $a\in\K$ and $\|a\|_\infty\leq 1$, and assume that its optimal value $\OPT$ is given to the algorithm as a public input, the same on neighboring instances. Let $\epsilon\in(0,1]$, $\delta\in(0,1/2]$, $\beta\in (0,1)$ and $\alpha\in(0,\OPT]$. Run Algorithm~\ref{alg:constraint_mwu} on $\K_{\OPT}=\{x\in\K:\Tr(x)=\OPT\}$ with width $\rho=\OPT$, accuracy $\alpha$, $T=\lceil 36\alpha^{-2}\OPT^2\log m\rceil$ rounds (with $m$ read as the public bound $\bar m$ of Remark~\ref{rem:public_rounds}, here and in the condition on $s$) and density parameter $s$, using the oracle $O(y)$ of Lemma~\ref{lem:private_oracle_eja} with $\alpha_0=\alpha/6$, privacy parameter $\epsilon'=\epsilon/\sqrt{8T\log(1/\delta)}$ and failure probability $\beta/T$. We also assume that every instance has at least $s$ constraints, so that the $1/s$-dense distributions used by Algorithm~\ref{alg:constraint_mwu} exist. There is an absolute constant $C_0$ such that, if
\begin{align*}
    s \geq & ~ C_0\cdot\frac{\OPT^2}{\epsilon\,\alpha^2}\cdot \log^{1/2}(1/\delta)\log^{1/2} m\cdot\Big(k_C\log\frac{2\OPT}{\alpha}+\log\frac{\log m}{\beta}\Big),
\end{align*}
then, for every fixed value of this input that is at least $1$ (Remark~\ref{rem:public_opt}), Algorithm~\ref{alg:constraint_mwu} is $(\epsilon,\delta)$-private with respect to high sensitivity constraint privacy, its output $x^*$ satisfies $x^*\in\K$ and $\Tr(x^*)=\OPT$, and with probability at least $1-\beta$ fewer than $s$ constraints have $\langle a_i, x^*\rangle<1-\alpha$.
\end{theorem}

\begin{proof}
We follow the proof of Theorem~\ref{thm:private_covering_sdp}. Recall that $\OPT\geq 1$. By self-duality of $\K$ and Lemma~\ref{lem:holder}, every $x\in\K_{\OPT}$ and every admissible constraint $a$ satisfy $0\leq\langle a,x\rangle\leq\|a\|_\infty\|x\|_1\leq\|x\|_1=\Tr(x)=\OPT$, so $|\langle a,x\rangle-1|\leq\OPT$. Hence $\rho=\OPT$ is a valid width, and $\K_{\OPT}$ contains an optimal solution. With $\alpha_0=\alpha/6$, privacy parameter $\epsilon'$ and failure probability $\beta/T$, the oracle of Lemma~\ref{lem:private_oracle_eja} is $\epsilon'$-private and $\OPT$-bounded, and its error is at most $\alpha/3$ provided
\begin{align*}
    \frac{6\OPT}{s\cdot \epsilon'}\Big(2k_C\log\frac{42\OPT}{\alpha}+\log\frac{T}{\beta}\Big) \leq & ~ \frac{\alpha}{6}.
\end{align*}
As in the proof of Theorem~\ref{thm:private_covering_sdp}, this condition is equivalent to
\begin{align*}
    s \geq & ~ \frac{216\sqrt{8}\,\OPT^2}{\alpha^2\epsilon}\cdot \sqrt{\log(1/\delta)\log m}\cdot \Big(2k_C\log\frac{42\OPT}{\alpha}+\log\frac{36\OPT^2\log m}{\alpha^2\beta}\Big).
\end{align*}
As $\log 21\leq 5\log 2\leq 5\log(2\OPT/\alpha)$, and as in that proof with $k_C\geq 1$ in place of $r\geq 1$, we get
\begin{align*}
    2k_C\log\frac{42\OPT}{\alpha} = & ~ 2k_C\Big(\log 21+\log\frac{2\OPT}{\alpha}\Big) \leq 12k_C\log\frac{2\OPT}{\alpha}, \\
    \log\frac{36\OPT^2\log m}{\alpha^2\beta} \leq & ~ 8k_C\log\frac{2\OPT}{\alpha}+\log\frac{\log m}{\beta}.
\end{align*}
So the last factor is at most $20\,(k_C\log(2\OPT/\alpha)+\log(\log m/\beta))$, and the condition in the theorem with $C_0=20\cdot 216\sqrt{8}$ suffices, up to the same rounding of $T$. The rest of the argument is the same as for Theorem~\ref{thm:private_covering_sdp}. It applies Lemma~\ref{lem:approximate_feasibility} to the constraints $\langle -a_i,x\rangle\leq -1$ and uses the convexity of $\K_{\OPT}$ and Theorem~\ref{thm:constraint_privacy_guarantee}.
\end{proof}

\section{Low Sensitivity SCPs}
\label{sec:low_sensitivity}

For low sensitivity SCPs, the divergence between neighboring inputs diminishes as the database size increases. We continue to study the feasibility program in the following form:
\begin{align*}
  & \mathrm{find} ~ x \in {\cal K} \\
  \mathrm{s.t.} ~ & \langle a_i, x\rangle\leq b_i, \forall i\in [m].
\end{align*}
We further without loss of generality normalize the constraints so that the feasible solutions are distributions over $\J$. The reduction is simple, if the optimal solution has trace $L$, then consider
\begin{align*}
  & \mathrm{find} ~ x \in {\cal K} \\
  \mathrm{s.t.} ~ & \langle a_i, x\rangle \leq b_i/L, \forall i\in [m]
\end{align*}
has a distribution solution. Once computing a solution $x^*$, we can scale back to obtain a solution for the unscaled program. The only downside is that if the constraints are only approximately satisfied: $\langle a_i,x^*\rangle \leq b_i/L + \alpha$ for all $i\in [m]$, by setting $\alpha = \alpha'/L$, then $\langle a_i, Lx^*\rangle \leq b_i + \alpha'$.

\subsection{Solving SCPs with Multiplicative Weights}

For convenience, we present a standard algorithm framework for multiplicative weights update. We note Algorithm~\ref{alg:mw} is quite different from the standard multiplicative weights update, which would take the form
\begin{align*}
    F^{t+1} \leftarrow \frac{F^t \circ \exp(-\eta\ell^t)}{\Tr(F^t\circ \exp(-\eta \ell^t))},
\end{align*}
however, this turns out to be insufficient since $\exp(A)\circ \exp(B)\neq \exp(A+B)$ in general unless $A$ and $B$ commute, see~\cite{clpv23} for more details. The regret bound of Algorithm~\ref{alg:mw} has been recorded in Theorem~\ref{thm:mw_regret}.
\begin{algorithm}[ht!] \caption{Multiplicative weights update algorithm.}\label{alg:mw}
  \begin{algorithmic}[1]
    \Procedure{MWU}{${\cal F}, \eta$}
    \State $F^1\gets \text{the uniform distribution on ${\cal F}$}$
    \For{$t=1\to T-1$}
        \State Receive loss $\ell^t\in \J$ (may depend on $F^1,\ldots,F^t$)
        \State $F^{t+1}\gets \frac{\exp(-\eta \sum_{i=1}^t \ell^i)}{\Tr(\exp(-\eta \sum_{i=1}^t \ell^i))}$
    \EndFor
    \State \Return $F^T$
    \EndProcedure
  \end{algorithmic}
\end{algorithm}

Our MWU algorithm for solving SCP is somewhat different from the variants introduced in~\cite{zvtl24}, as it relies on a dual oracle that returns an approximately most violated constraint. This should be treated as a generalization of the dual oracle for LP introduced in~\cite{hrru14}.
\begin{definition}[Dual oracle]
For $\gamma>0$, an $(\alpha,\gamma)$-dual oracle, given
$A=\{a_1,\ldots,a_m\}\in\mathcal{J}^m$, $b\in\mathbb{R}^m$, and $x\in\mathcal{K}$ as input,
returns an index $i\in[m]$ such that
\begin{align*}
  \langle a_i, x\rangle - b_i \ge \max_{j\in[m]} (\langle a_j, x\rangle - b_j) - \alpha,
\end{align*}
with probability at least $1-\gamma$.
\end{definition}

The full algorithm is given in~\ref{alg:primal_lp_nonpriv}.
 
\begin{algorithm}[ht!]\caption{Solving SCP feasibility via MWU over primal variables.}
  \label{alg:primal_lp_nonpriv}
  \begin{algorithmic}[1]
  \Procedure{SCPMWUPrimal}{$A \in {\cal J}^m$, $b \in \R^m$}
    \State $x^1\gets e/r$
    \State Let $\rho =
      \max_{i\in [m]} \|a_i\|_\infty$ be the width of the SCP, $\alpha > 0$ be the
      desired accuracy
    \State Let $\textsc{Oracle}$ be a $(\alpha, \gamma)$-dual
      oracle, 
      \State $\eta \gets \frac{\alpha}{4\rho}$, $  T \gets \frac{16 \rho^2 \log r}{\alpha^2}$ 
    \For{ $t = 1 \to T$}
        \State  $p^t \gets \textsc{Oracle}(A, b, x^t)$
        \State $\ell^t \gets \frac{1}{\rho} a_{p^t}$
        \State $x^{t+1}\gets \frac{\exp(-\sum_{i=1}^t\eta \ell^i)}{\Tr(\exp(-\sum_{i=1}^t\eta \ell^i))}$.
    \EndFor
    \State{\Return $\ov{x} = (1/T) \sum_{t = 1}^T x^t$}
    \EndProcedure
  \end{algorithmic}
  
\end{algorithm}

Our MWU framework is based on SCP variant introduced in~\cite{clpv23,zvtl24}, hence its proof differs significantly from~\cite{pst95} as the update is different. Moreover, the SCP solver due to~\cite{zvtl24} utilizes \emph{primal oracle}. Hence, we prove a convergence theorem for the dual oracle.

\begin{theorem}[Formal version of Theorem~\ref{thm:non_private_scp_informal}]
\label{thm:non_private_scp}
Suppose the SCP admits a distributional feasible solution $x \in \mathcal{K}$ with $\Tr(x) = 1$. 
Then Algorithm~\ref{alg:primal_lp_nonpriv}, when equipped with an $(\alpha/2,\gamma)$-dual oracle, 
returns $x \in \mathcal{K}$ such that
\begin{align*}
   \langle a_i, x \rangle \leq b_i + \alpha \qquad \forall i \in [m],
\end{align*}
with probability at least $1 - T\gamma$.
\end{theorem}

\begin{proof}
Our proof will be applying Theorem~\ref{thm:mw_regret}, to do so, we first need to verify that the sequence of loss satisfy that $\|\ell_t\|_\infty\leq 1$. As we set $\rho$ to be the max infinity norm over all constraints, it's easy to see that $\|\ell^t\|_\infty\leq 1$. Next, we assume the oracle is exact, i.e., it returns $\max_{j\in [m]}\langle a_j,x\rangle-b_j$. Then we show how to generalize the argument to approximate oracles. By Theorem~\ref{thm:mw_regret}, we have that
\begin{align}\label{eq:mw_regret}
    \sum_{t=1}^T \langle \ell^t, x^t\rangle \leq & ~ \sum_{t=1}^T \langle \ell^t, y\rangle+\eta T+\frac{\ln r}{\eta},
\end{align}
for any probability distribution $y\in {\cal K}$. Set $r^t:=\langle a_{p^t},x^t\rangle-b_{p^t}$, we can explicitly compute the LHS:
\begin{align*}
    \sum_{t=1}^T \langle \ell^t,x^t\rangle= & ~ \frac{1}{\rho}\sum_{t=1}^T \langle a_{p^t},x^t\rangle \\
    = & ~ \frac{1}{\rho}\sum_{t=1}^T (b_{p^t}+r^t),
\end{align*}
meanwhile since we can choose any distribution $y\in {\cal K}$, we choose a feasible $y$ so that for any $i\in [m]$, $\langle a_i,y\rangle\leq b_i$, thus
\begin{align*}
    \sum_{t=1}^T \langle \ell^t,y\rangle = & ~ \frac{1}{\rho}\sum_{t=1}^T \langle a_{p^t},y\rangle \\
    \leq & ~ \frac{1}{\rho}\sum_{t=1} b_{p^t},
\end{align*}
combining these inequalities yields
\begin{align*}
    \frac{1}{\rho}\sum_{t=1}^T r^t \leq & ~ \eta T+\frac{\ln r}{\eta},
\end{align*}
multiplying both sides by $\frac{\rho}{T}$ gives
\begin{align*}
    \frac{1}{T}\sum_{t=1}^T r_t \leq & ~ \eta \rho+\frac{\rho\ln r}{\eta T},
\end{align*}
setting $\eta=\frac{\epsilon}{2\rho}$ and $T=\frac{4\rho^2\ln r}{\epsilon^2 }$, the above bound simplifies to
\begin{align*}
    \frac{1}{T}\sum_{t=1} r_t \leq & ~ \epsilon,
\end{align*} 
to obtain a final bound for any constraint, fix $i\in [m]$, we note that for any $t\in [T]$, it must be the case that $\langle a_i,x^t\rangle-b_i\leq \langle a_{p^t},x^t\rangle-b_{p^t}$, and if we average it gives
\begin{align*}
    \frac{1}{T}\sum_{t=1}^T \langle a_i,x^t\rangle-b_i = & ~ \langle a_i,\ov x\rangle-b_i \\
    \leq & ~ \frac{1}{T}\sum_{t=1}^T r_t \\
    \leq & ~ \epsilon,
\end{align*}
as desired. Now recall in our algorithm we are using an approximate dual oracle, and we slightly modify the analysis by setting $\epsilon=\alpha/2$, and note that the only place we are using the most violated constraint is to prove all constraints, where using an $(\alpha/2,\gamma)$-dual oracle would blow up the final error by a factor of $\alpha/2$. Together with the choice of $\epsilon$, we conclude that for any $i\in [m]$, $\langle a_i,\ov x\rangle-b_i\leq \alpha$ holds with probability at least $1-T\gamma$ (via a union bound).
\end{proof}

\subsection{Scalar Private SCPs} \label{sec:qr}

We first consider the ``simpler'' case of scalar private SCPs. This is when the objective and constraints are public data, and only the scalar vector $b$ is private. A private database $D$ is mapped to a tuple $(A(D), c(D), b(D))$. 
For neighboring databases $D,D'$, the mappings satisfy
\begin{align*}
   c(D) = & ~ c(D'), \\
   A(D) = & ~  A(D'),
\end{align*}
i.e., both the objective vector $c$ and the constraint matrix $A$ are independent of the data. Moreover, the right-hand side vectors differ only by a bounded amount:
\begin{align*}
   \| b(D) - b(D') \|_\infty \leq & ~  \Delta_\infty .
\end{align*}
 Note that since $b\in \R^m$, here the infinity norm is the standard vector $\ell_\infty$ norm. We again focus on the feasibility SCP. Formally:

\begin{definition}
Given a vector $b \in \mathbb{R}^m$, and a constraint set $A \in \J^m$, 
a randomized mechanism $\mathcal{M}$ that outputs a vector in $\J$ is $(\epsilon, \delta)$-low sensitivity scalar private if for any $b, b'$ with $\|b-b'\|_\infty \leq \Delta_\infty$,
\begin{align*}
    \Pr[\mathcal{M}(b, A) \in S] \leq e^\epsilon \Pr[\mathcal{M}( b', A) \in S] + \delta
\end{align*}
for any subset $S \subseteq \J$.
\end{definition}

To achieve privacy, we will implement a private dual oracle via exponential mechanism, henceforce privatize Algorithm~\ref{alg:primal_lp_nonpriv}.

\begin{lemma} \label{lem:qp_dual}
Let $\epsilon > 0$ and $\gamma \in (0,1)$. Suppose that on neighboring instances 
the vector $b$ changes by at most $\Delta_\infty$ in $\ell_\infty$ norm. 
Then, the $\epsilon$-private exponential mechanism with quality score
\begin{align*}
   Q(i,b) = \langle a_i, x \rangle - b_i
\end{align*}
is an $(\alpha,\gamma)$-dual oracle, where
\begin{align*}
   \alpha = \frac{2 \Delta_\infty}{\epsilon}\log(\frac{m}{\gamma}).
\end{align*}
\end{lemma}

\begin{proof}
    Note that our oracle will be following the procedure:
    \begin{itemize}
        \item Given $A, x, b$, computing $Q(i, b)$ for all $i\in [m]$;
        \item Sample constraint $i$ with probability $\exp(\frac{\epsilon}{2\Delta_\infty}\cdot Q(i,b))$.
    \end{itemize}
    This procedure is automatically $\epsilon$-private following the definition of exponential mechanism (Def.~\ref{def:exponential_mechanism}), so it remains to prove the procedure indeed implements an $(\alpha,\gamma)$-dual oracle. By Lemma~\ref{lem:acc_gua_exp_mec}, this is true if $\alpha=\frac{2\Delta_\infty}{\epsilon}\cdot \log(\frac{m}{\gamma})$, which completes the proof.
\end{proof}

We are now in the position to provide privacy and accuracy guarantees for scalar private SCPs.

\begin{theorem} \label{thm:qp_priv}
Suppose the SCP admits a feasible $x\in\K$ with $\Tr(x)=1$. Let $\epsilon \in (0,1]$, $\delta\in(0,1/2]$, $\beta \in (0,1)$, and define 
$\rho = \max_{i \in [m]} \|a_i\|_\infty$. 
Then Algorithm~\ref{alg:primal_lp_nonpriv} is 
$(\epsilon,\delta)$-low sensitivity scalar private with sensitivity $\Delta_\infty$, 
and with probability at least $1-\beta$ outputs $x^* \in \mathcal{K}$ such that $\langle a_i,x^*\rangle\leq b_i+\alpha$ for all $i\in [m]$, where
    \begin{align*}
        \alpha = & ~\wt O\left(\frac{\rho^{1/2}\Delta_\infty^{1/2}}{\epsilon^{1/2}}\log^{1/4} r\log^{1/4}(1/\delta)\log^{1/2}(1/\beta)\log^{1/2}m\right).
    \end{align*}
\end{theorem}

\begin{proof}
    Let $\epsilon'=\frac{\epsilon}{\sqrt{8T\log(1/\delta)}}$, we use exponential mechanism to implement an $\epsilon'$-private dual oracle, and by adaptive composition (Lemma~\ref{lem:composition_lemma}), this gives a final $(\epsilon,\delta)$-private algorithm. For success probability, we set $\gamma=\beta/T$. Theorem~\ref{thm:non_private_scp} needs an $(\alpha/2,\gamma)$-dual oracle, and by Lemma~\ref{lem:qp_dual} the exponential mechanism is one when
    \begin{align*}
        \frac{\alpha}{2} = & ~ \frac{2\Delta_\infty}{\epsilon'} \cdot \log(\frac{mT}{\beta}) \\
        = & ~ \frac{2\Delta_\infty \sqrt{8T\log(1/\delta)}}{\epsilon}\cdot \log(\frac{mT}{\beta}) \\
        = & ~ \frac{16\sqrt{2}\,\Delta_\infty \rho\log^{1/2}r\log^{1/2}(1/\delta)}{\alpha \epsilon}\cdot \log(\frac{16m\rho^2 \log r}{\alpha^2 \beta}),
    \end{align*}
    solving for $\alpha$ gives
    \begin{align*}
        \alpha = & ~ \wt O\left(\frac{\rho^{1/2}\Delta_\infty^{1/2}}{\epsilon^{1/2}}\log^{1/4} r\log^{1/4}(1/\delta)\log^{1/2}(1/\beta)\log^{1/2}m\right),
    \end{align*}
    as desired. This completes the proof.
\end{proof}

\subsection{Low Sensitivity Constraint Private SCPs} \label{sec:mp}
 
Given the feasibility program
\begin{align*}
  &\mathrm{find} ~ x\in {\cal K} \\
  \mathrm{s.t.} ~ & \langle a_i,x\rangle\leq b_i, \forall i\in [m],
\end{align*}
and a neighboring instance would perturb some constraints by small amounts. More specifically, for two set of constraints $A, A'\in \J^m$, we define the global infinity norm as
\begin{align*}
    \|A-A'\|_\infty = & ~ \max_{i\in [m]} \|a_i-a'_i\|_\infty,
\end{align*}
and we assume for two neighboring instances, their distance in global infinity norm is at most $\Delta_\infty$. 

Similarly to the scalar private setting, a private database $D$ is mapped to a tuple $(A(D), c(D), b(D))$. 
For every pair of neighboring databases $D,D'$, the mappings satisfy
\begin{align*}
   c(D) = & ~ c(D'), \\
   b(D) = & ~  b(D'),
\end{align*}
i.e., both the objective vector $c$ and the scalar $b$ are independent of the data. Moreover, the constraints differ by a bounded amount
\begin{align*}
   \| A(D) - A(D') \|_\infty \leq & ~  \Delta_\infty .
\end{align*}
Formally:

\begin{definition}
Given a vector $b \in \mathbb{R}^m$, and a constraint set $A \in \J^m$, 
a randomized mechanism $\mathcal{M}$ that outputs a vector in $\J$ is $(\epsilon, \delta)$-low sensitivity constraint private if for any $A, A'$ such that $\|A-A'\|_\infty \leq \Delta_\infty$,
\begin{align*}
    \Pr[ {\cal M}(b, A) \in S ] \leq e^\epsilon \Pr[ {\cal M}(b, A') \in S] + \delta 
\end{align*}
for any subset $S \subseteq \J$.
\end{definition}

We without loss of generality normalize the constraints so that the spectrum of each $a_i$ lies in $[-1, 1]$. Our algorithm will again be implementing the dual oracle with exponential mechanism. However, as the oracle returns a constraint and it will be used to compute the loss, we have to add one more layer of privacy via the generic Gaussian mechanism.

\begin{theorem}[Privacy guarantee of Theorem~\ref{thm:rowsens_informal}]
\label{thm:rowsens_priv}
Let $\epsilon\in(0,1]$ and $\delta\in(0,1/2]$, let $\epsilon'$ and $\Delta_\infty$ be as defined in 
Algorithm~\ref{alg:primal_lp_rowsens}, and suppose the algorithm 
employs an $\epsilon'$-private dual oracle. 
Then Algorithm~\ref{alg:primal_lp_rowsens} is 
$(\epsilon,\delta)$-low sensitivity constraint private 
with sensitivity $\Delta_\infty$.
\end{theorem}

\begin{proof}
We note the privacy comes from two sources: the Gaussian mechanism and the oracle operation. For Gaussian mechanism, the noise is added to $a_{p^t}$, whose $\ell_2$ sensitivity is at most $\sqrt{r}\Delta_\infty$, so by Lemma~\ref{lem:gaussian_mech_eja} (with $\epsilon'<1$ and $\sigma$ as in Algorithm~\ref{alg:primal_lp_rowsens}) each operation is $(\epsilon', \delta/(4T))$-private, and each oracle is $\epsilon'$-private. By an adaptive composition (Lemma~\ref{lem:composition_lemma}) over $2T$ operations with slack $\delta/2$, using $\epsilon'=\epsilon/\sqrt{16T\log(2/\delta)}$, we see that the algorithm is $(\epsilon,\delta/2+2T\cdot\delta/(4T))=(\epsilon,\delta)$-private.
\end{proof}

We prove that the exponential mechanism is a private dual oracle under this setting of neighboring.

\begin{lemma} \label{lem:rp_dual}

Let $\epsilon > 0$ and $\gamma \in (0,1)$. Suppose that on neighboring instances 
the constraint set $A$ changes by at most $\Delta_\infty$ in $\ell_\infty$ norm. Let $x\in \K$ be any distributional element. Then, the $\epsilon$-private exponential mechanism with quality score
  \begin{align*}
    Q(i, A) = \langle a_i,x\rangle - b_i
  \end{align*}
  is an $(\alpha, \gamma)$-dual oracle, for
  \begin{align*}
    \alpha = \frac{2 \Delta_\infty}{\epsilon} \cdot \log  ( \frac{m}{\gamma}
     ).
  \end{align*}
\end{lemma}
\begin{proof}
Since $x$ is a distribution, the quality score $Q$ changes by at most 
$\Delta_\infty$ on neighboring inputs, and hence is $\Delta_\infty$-sensitive. 
The claimed accuracy bound then follows directly from the accuracy guarantee 
of the exponential mechanism (Lemma~\ref{lem:acc_gua_exp_mec}).
\end{proof}

To prove the convergence, we need to slightly change the original argument, as we have to perturb each constraint by a Gaussian noise element. We give a customized proof based on standard regret bound.

\begin{theorem}[Utility guarantee of Theorem~\ref{thm:rowsens_informal}]
\label{thm:rowsens_acc}
Let $A \in \mathcal{J}^m$ satisfy $\lambda(a_i) \subseteq [-1,1]$ for all $i \in [m]$, 
and let $b \in \mathbb{R}^m$. 
Fix $\beta\in(0,1)$, $\epsilon\in(0,1]$ and $\delta \in (0,1/2]$, and suppose the SCP admits a feasible $y\in\mathcal{K}$ with $\Tr(y)=1$. 
Let $r$ denote the rank of $\mathcal{J}$ and $k$ its dimension. 
Then, with probability at least $1-\beta$, 
Algorithm~\ref{alg:primal_lp_rowsens}, when run with the exponential mechanism 
as a dual oracle (Lemma~\ref{lem:rp_dual}), 
returns a distribution $x^*$ such that
\begin{align*}
   \langle a_i, x^* \rangle \leq b_i + \alpha 
   \qquad \forall i \in [m],
\end{align*}
where
\begin{align*}
   \alpha = & ~ \widetilde{O}\left(
     \frac{\Delta_\infty^{1/2} r^{1/4} k^{1/4}}{\epsilon^{1/2}}
     \cdot \mathrm{poly}\log(r, m, 1/\beta, 1/\delta)
   \right).
\end{align*}
\end{theorem}

\begin{proof}
  Let $\epsilon'$ be as in Theorem~\ref{thm:rowsens_priv}, 
$T$ be the number of iterations in Algorithm~\ref{alg:primal_lp_rowsens}, 
and set $\gamma = \beta/(2T)$. 
By Lemma~\ref{lem:rp_dual}, with probability at least $1-\gamma$, 
the oracle returns $p^t$ such that for all constraints $i \in [m]$,
\begin{align}\label{eq:most_violated_constraint}
     ( \langle a_i, x^t\rangle - b_i )  -  (
      \langle a_{p^t}, x^t\rangle - b_{p^t}  ) \leq \frac{2
      \Delta_\infty}{\epsilon'} \cdot \log  ( \frac{m}{\gamma}  )
  \end{align}
    Here $p^t$ is chosen as the constraint that is nearly the most violated, 
up to an additive factor of $\alpha$. 
Define the vanilla loss $\ell^t = a_{p^t}$. 
Note that the left-hand side of \eqref{eq:most_violated_constraint} 
is exactly $(\langle a_i,x^t\rangle-b_i)-(\langle \ell^t,x^t\rangle-b_{p^t})$. Applying a union bound over the $T$ oracle calls, 
inequality~\eqref{eq:most_violated_constraint} 
holds simultaneously for all $t \in [T]$ with probability at least $1-\beta/2$. 
We henceforth condition on this event.

We next need to bound the norm of the noise element $z$. Again, we utilize the isometric between $\J$ and $\R^k$ in $\ell_2$ norm: recall that $\phi(z^t)$ is a Gaussian vector sampled from $\N(0, \sigma^2 I_k)$, by Lemma~\ref{lem:chi_square_tail}, with probability at least $1-\beta/(2T)$,
  \begin{align}\label{eq:gaussian_norm}
      \|z^t\|_2 \leq & ~ \sigma \cdot (\sqrt k+\sqrt {2\log(\frac{2T}{\beta})})\notag \\
      = & ~ \frac{\Delta_\infty(\sqrt{2rk\log(5T/\delta)}+2\sqrt{r\log(5T/\delta)\log(1/\gamma)})}{\epsilon'}
  \end{align}
  union bound over $T$ noise elements, it succeeds with probability at least $1-\beta/2$, condition on this event happen.
  
  We proceed by assuming the norm of Gaussian noise and the error of the oracle are small, in particular,
  \begin{align}\label{eq:small_noise}
        \frac{\Delta_\infty(\sqrt{2rk\log(5T/\delta)}+2\sqrt{r\log(5T/\delta)\log(1/\gamma)})}{\epsilon'} \leq & ~ \frac{\alpha}{6},
  \end{align}
  and
  \begin{align}\label{eq:small_oracle_error}
        \frac{2\Delta_\infty}{\epsilon'}\log\Big(\frac{m}{\gamma}\Big) \leq & ~ \frac{\alpha}{6},
  \end{align}
  since $\alpha<1$, Eq.~\eqref{eq:small_noise} implies that the norm of the noise is at most $1/6$. Hence, 
  \begin{align*}
      \|\wh \ell^t\|_\infty\leq & ~ \frac{\|\ell^t \|_\infty+\|z^t\|_\infty}{2} \\
      \leq & ~  \frac{\|\ell^t \|_\infty+\|z^t\|_2}{2} \\
      \leq & ~ 1.
  \end{align*}
   We can then apply Theorem~\ref{thm:mw_regret}: pick $y\in {\cal K}$ be a feasible distribution element, then
  \begin{align*}
      \frac{1}{T}\sum_{t=1}^T \langle \wh \ell^t,x^t\rangle \leq & ~ \frac{1}{T}\sum_{t=1}^T \langle \wh \ell^t, y\rangle+\eta+\frac{\log r}{\eta T},
  \end{align*}
set $r_t:=\frac{\langle a_{p^t},x^t\rangle-b_{p^t}}{2}+\frac{\langle z^t, x^t\rangle}{2}$, we compute
\begin{align*}
    \frac{1}{T}\sum_{t=1}^T r_t = & ~ \frac{1}{T}\sum_{t=1}^T \langle \wh \ell^t,x^t\rangle-\frac{b_{p^t}}{2} \\
    \leq & ~ \frac{1}{T}\sum_{t=1}^T \langle \wh\ell^t,y\rangle-\frac{b_{p^t}}{2}+\eta+\frac{\log r}{\eta T}.
\end{align*}
Fix any constraint $a_i$, since the error of the exponential mechanism is at most $\frac{2\Delta_\infty}{\epsilon'}\log(\frac{m}{\gamma})\leq\frac{\alpha}{6}$ by Eq.~\eqref{eq:small_oracle_error} and $x^t$ is a distribution, it must be the case that
\begin{align*}
    \frac{1}{T}\sum_{t=1}^T \frac{\langle a_i,x^t\rangle-b_i}{2}+\frac{\langle z^t,x^t\rangle}{2} \leq & ~ \frac{1}{T}\sum_{t=1}^T r_t+\frac{\alpha}{6} \\
    \leq & ~ \frac{1}{T}\sum_{t=1}^T \langle \wh\ell^t,y\rangle-\frac{b_{p^t}}{2}+\eta+\frac{\log r}{\eta T}+\frac{\alpha}{6},
\end{align*}
where the first step is by Eq.~\eqref{eq:most_violated_constraint}. Recall that we also have shown that the norm of Gaussian noise is small and we could bound the term introduced by it:
\begin{align*}
    |\langle z^t,x^t\rangle | \leq & ~ \|x^t\|_1\cdot \|z^t \|_\infty  \\
    \leq & ~ \alpha/6,
\end{align*}
where the first step is by H{\"o}lder's inequality (Lemma~\ref{lem:holder}) and the second step is by Eq.~\eqref{eq:small_noise}. Thus, we have
\begin{align*}
    \frac{1}{T}\sum_{t=1}^T \langle a_i,x^t\rangle-b_i \leq & ~ \frac{1}{T}\sum_{t=1}^T 2\langle\wh \ell^t,y\rangle-b_{p^t}+\frac{\alpha}{2}+2\eta+\frac{2\log r}{\eta T}.
\end{align*}
It remains to examine $2\langle\wh \ell^t,y\rangle-b_{p^t}$, since $y$ is feasible, it must be that for any $i$, $\langle a_i,y\rangle-b_i\leq 0$, therefore
\begin{align*}
    2\langle\wh \ell^t,y\rangle-b_{p^t} = & ~ \langle a_{p^t},y\rangle+\langle z^t, y\rangle-b_{p_t} \\
    \leq & ~ \langle z^t, y\rangle \\
    \leq & ~ \|y\|_1\cdot \|z^t\|_\infty \\
    \leq & ~ \alpha/6,
\end{align*}
put things together, we have
\begin{align*}
    \frac{1}{T}\sum_{t=1}^T \langle a_i,x^t\rangle-b_i = & ~ \langle a_i,\ov x\rangle-b_i \\
    \leq & ~ \frac{2\alpha}{3}+2\eta +\frac{2\log r}{\eta T} \\
    \leq & ~ \alpha,
\end{align*}
where the last step is by the choice of $\eta$ and $T$. It remains to provide a value of $\alpha$ that satisfies Eq.~\eqref{eq:small_noise} and Eq.~\eqref{eq:small_oracle_error}. Since $\frac{1}{\epsilon'}=\frac{4\sqrt{T\log(2/\delta)}}{\epsilon}$, $\log(2/\delta)\leq\log(5T/\delta)$ and $\sqrt{T}=\frac{12\sqrt{\log r}}{\alpha}$, Eq.~\eqref{eq:small_noise} holds if

  \begin{align*}
    & ~ \frac{6\Delta_\infty(\sqrt{2rk\log(5T/\delta)}+2\sqrt{r\log(5T/\delta)\log(1/\gamma)})}{\epsilon'} \\
    \leq & ~ \frac{48\Delta_\infty \sqrt{rT}\cdot (\sqrt{k}+\sqrt{\log(1/\gamma)})\cdot \log(5T/\delta)}{\epsilon} \\
    = & ~ \frac{576\Delta_\infty  \sqrt{r\log r}\cdot (\sqrt{k}+\sqrt{\log(2T/\beta)})\cdot \log(5T/\delta)}{\alpha\epsilon} \leq \alpha,
  \end{align*}
  and, using $\sqrt{k}+\sqrt{L}\leq 2\sqrt{kL}$ for $k\geq 1$ and $L=\log(2T/\beta)\geq 1$, this holds if the first inequality below holds. In the same way, Eq.~\eqref{eq:small_oracle_error} is equivalent to the second inequality:
  \begin{align*}
      \alpha \geq & ~ \frac{34\Delta_\infty^{1/2}r^{1/4}(\log r)^{1/4}k^{1/4}}{\epsilon^{1/2}}\cdot \log^{1/4}(2T/\beta)\cdot \log^{1/2}(5T/\delta), \\
      \alpha \geq & ~ \frac{24\Delta_\infty^{1/2}(\log r)^{1/4}}{\epsilon^{1/2}}\cdot \log^{1/4}(2/\delta)\cdot \log^{1/2}(2mT/\beta).
  \end{align*}
  As $T=144\log r/\alpha^2$ enters only through logarithms, both hold for some
  \begin{align*}
      \alpha = & ~ \widetilde O\left(\frac{\Delta_\infty^{1/2} r^{1/4} k^{1/4}}{\epsilon^{1/2}}\cdot\mathrm{poly}\log(r,m,1/\beta,1/\delta)\right).
  \end{align*}
 This completes the proof.
\end{proof}

\begin{remark}
Our algorithm differs significantly from the row private LP algorithm of~\cite{hrru14}, where they ensure the constraints are private by injecting Laplace noises to each entry. This comes from the fact that for $\J=\R^r$, the only Jordan frame is the standard basis, hence to develop a differential private mechanism, it is enough to perturb the eigenvalues. For other $\J$, this is no longer the case. Consider $\J$ to be the set of all $r\times r$ real symmetric matrices, then a Jordan frame can be formed by taking any set of orthonormal vectors $u_1,\ldots,u_r\in \R^r$ and computing $u_1u_1^\top,\ldots,u_ru_r^\top$. This large degree of freedom means that private mechanism must also perturb the Jordan frame to ensure the basis information is preserved. We achieve so by implicitly resorting to the isomorphism between $\J$ and $\R^k$, as such a random Gaussian element would possess both random eigenvalues and Jordan frame. Note that we choose Gaussian mechanism instead of Laplace mechanism, as we will only provide either $\ell_2$ or $\ell_1$ norm guarantee instead of $\ell_\infty$, and $\ell_2$ norm only distorts $\ell_\infty$ by a factor of $\sqrt k$ instead of $k$ that is given by the $\ell_1$ norm.
\end{remark}

\subsection{Objective Private SCPs} \label{sec:op}

Finally, we consider objective private SCP, where instead of solving a feasibility problem, we try to solve the optimization version. In contrast to~\cite{hrru14}, we again consider two neighboring objectives differ in $\ell_\infty$ norm. 
\begin{align*}
  \max_{x\in \J} & ~ \langle c,x\rangle \\
  \text{s.t.} & ~ \langle a_i, x\rangle\leq b_i, \forall i\in [m] \\
  & ~ x\in \K.
\end{align*}
Given two neighboring databases $D, D'$, we have $A(D)=A(D'), b(D)=b(D')$ and $\|c(D)-c(D')\|_\infty\leq \Delta_\infty$. Formally,

\begin{definition}
Given a vector $b \in \mathbb{R}^m$, $c\in \K$ and a constraint set $A \in \J^m$, 
a randomized mechanism $\mathcal{M}$ that outputs a vector in $\J$ is $(\epsilon, \delta)$-low sensitivity objective private if for any $c, c'$ such that $\|c-c'\|_\infty \leq \Delta_\infty$,
\begin{align*}
    \Pr[ {\cal M}(c, b, A) \in S ] \leq e^\epsilon \Pr[ {\cal M}(c', b, A) \in S] + \delta 
\end{align*}
for any subset $S \subseteq \J$.
\end{definition}

Next, we present a simple algorithm, based similarly on the Gaussian mechanism introduced in Lemma~\ref{lem:gaussian_mech_eja}, but applied to the objective element. We further assume a somewhat unusual but necessary condition on the SCP: there exists an optimal solution with unit $\ell_2$ norm. 
\begin{theorem}
    Let the objective private SCP has optimal value ${\rm OPT}$ and the optimal solution has unit $\ell_2$ norm. Suppose $\dim(\J)=k$ and $\phi:\J\rightarrow \R^k$ is an isometry between $\J$ and $\R^k$. Consider the following mechanism:
    \begin{itemize}
        \item Set $\sigma=\frac{\Delta_\infty\sqrt{2r\log(1.25/\delta)}}{\epsilon}$, where $\epsilon,\delta\in(0,1)$;
        \item Generate a Gaussian noise vector $\nu\sim \N(0, \sigma^2I_k)$;
        \item Set $z=\phi^{-1}(\nu)$.
    \end{itemize}
    Then, let $\wt c:=c+z$, and the perturbed SCP
    \begin{align*}
        \max_{x\in \J} & ~ \langle \wt c, x\rangle \\
        \text{s.t.} & ~ \langle a_i, x\rangle\leq b_i, \forall i\in [m] \\
        & ~ \|x\|_2 = 1 \\
        & ~ x\in \K
    \end{align*}
is released. Then, the algorithm is $(\epsilon,\delta)$-low sensitivity objective private with sensitivity $\Delta_\infty$. With probability $1-\beta$, solving the perturbed SCP non-privately produces $x^*$ such that $\langle a_i,x^*\rangle\leq b_i, \forall i\in [m]$ and $\langle c,x^*\rangle\geq {\rm OPT}-\alpha$, where
\begin{align*}
    \alpha = & ~ \frac{4\Delta_\infty\sqrt{r\log(1.25/\delta)}(\sqrt k+\sqrt{\log(1/\beta)})}{\epsilon}.
\end{align*}
\end{theorem}

\begin{proof}
  Recall that $\|c\|_2\leq \sqrt r\cdot \|c\|_\infty$, by Corollary~\ref{cor:linfty_eja}, indeed the mechanism is $(\epsilon,\delta)$-private. For accuracy, since $\phi$ is an isometry, $\|z\|_2=\|\nu\|_2$, so by Lemma~\ref{lem:chi_square_tail}, with probability at least $1-\beta$,
  \begin{align*}
      \|z\|_2 \leq & ~ \sigma\cdot (\sqrt k+\sqrt{2\log(1/\beta)}) \\
      \leq & ~ \frac{2\Delta_\infty\sqrt{r\log(1.25/\delta)}(\sqrt k+\sqrt{\log(1/\beta)})}{\epsilon} = \frac{\alpha}{2}.
  \end{align*}
  Let $\wt x^*$ be the optimal solution to the perturbed SCP, and let $x^*$ be the optimal solution to the original SCP. Note that if
  \begin{align*}
      \langle c,\wt x^*\rangle< & ~ {\rm OPT}-\alpha,
  \end{align*}
  then
  \begin{align*}
      \langle \wt c, \wt x^*\rangle = & ~ \langle c+z, \wt x^*\rangle \\
      < & ~ {\rm OPT}-\alpha+\langle z,\wt x^*\rangle \\
      \leq & ~ {\rm OPT}-\alpha+\|z\|_2\cdot \|\wt x^*\|_2 \\
      \leq & ~ {\rm OPT}-\alpha/2,
  \end{align*}
  where we use Cauchy-Schwarz inequality and $\|z\|_2\leq\alpha/2$. Meanwhile, note that
  \begin{align*}
      \langle \wt c,x^*\rangle = & ~ {\rm OPT}+\langle z,x^*\rangle \\
      \geq & ~ {\rm OPT}-\|z\|_2\cdot \|x^*\|_2 \\
      \geq & ~ {\rm OPT}-\alpha/2,
  \end{align*}
  contradicts the definition of $\wt x^*$. Thus, it must be the case that $\langle c,\wt x^*\rangle\geq {\rm OPT}-\alpha$, and $\wt x^*$ is feasible since it's a feasible solution to the perturbed SCP, and the perturbation does not change feasibility. This completes the proof.
\end{proof}

\begin{remark}
In~\cite{hrru14}, they impose a much standard assumption that $\|x\|_1=1$. This again comes from the fact that for $\J=\R^k$, the isomorphism is just the identity map, and in addition to the isometry in $\ell_2$ norm $\|\phi(x)\|_2=\|x\|_2$, all norms are preserved. This is particularly important when applying H{\"o}lder's inequality: when $\|x\|_1=1$, we could upper or lower bound the inner product by the $\ell_\infty$ norm of the noise $z$. This is no longer true for $\J\neq \R^k$, as $\phi(\cdot)$ only preserves the $\ell_2$ norm, translating between different norms would incur blowup or shrinkage factors dependent on $r$. Of course, imposing an $\ell_2$ norm constraint makes the constraint set no longer an affine subspace, which would require solvers that could handle quadratic constraints.
\end{remark}


\ifdefined\isarxiv
\else
\section{Impact Statement}
\label{sec:impact}

Our work concerns the privacy of symmetric cone programming, we believe its development would help protect the data and model privacy for various machine learning tasks. At its current stage, our work is theoretical, so it does not lead to any direct negative societal consequences.

\fi

\end{document}